\documentclass[12pt]{amsart}
\usepackage{amssymb}
\usepackage{amsfonts}
\usepackage{color}
\newtheorem{thm}{Theorem}[section]

\newtheorem{lem}[thm]{Lemma}
\newtheorem{prop}[thm]{Proposition}
\theoremstyle{definition}
\newtheorem{defn}[thm]{Definition}
\theoremstyle{remark}
\newtheorem{rem}[thm]{Remark}
\newcommand{\R}{\mathbb R}
\newcommand{\C}{\mathbb C}

\newcommand{\X}{\mathfrak X}
\newcommand{\B}{\mathfrak B}
\newcommand{\K}{\mathfrak K}
\newcommand{\bb}{\mathfrak b}
\begin{document}

\title[WARPED PRODUCT K\"AHLER MANIFOLDS...]
{WARPED PRODUCT K\"AHLER MANIFOLDS AND BOCHNER-K\"AHLER METRICS}%
\author{G. Ganchev and V. Mihova}%
\address{Bulgarian Academy of Sciences, Institute of Mathematics and
Informatics, Acad. G. Bonchev Str. bl. 8, 1113 Sofia, Bulgaria}%
\email{ganchev@math.bas.bg}%
\address{Faculty of Mathematics and Informatics, University of Sofia,
J. Bouchier Str. 5, (1164) Sofia, Bulgaria}
\email{mihova@fmi.uni-sofia.bg}
\subjclass{Primary 53B35, Secondary 53C25}%
\keywords{$\alpha$-Sasakian manifolds, warped product K\"ahler
manifolds, Bochner-K\"ahler manifolds with special scalar distribution,
warped product Bochner-K\"ahler metrics.}%

\begin{abstract}
\vskip 2mm
Using as an underlying manifold an alpha-Sasakian manifold we introduce
warped product Kaehler manifolds. We prove that if the underlying
manifold is an alpha-Sasakian space form, then the corresponding
Kaehler manifold is of quasi-constant holomorphic sectional curvatures
with special distribution. Conversely, we prove that any Kaehler
manifold of quasi-constant holomorphic sectional curvatures with
special distribution locally has the structure of a warped product
Kaehler manifold whose base is an alpha-Sasakian space
form. Considering the scalar distribution generated by the scalar
curvature of a Kaehler manifold, we give a new approach to the local
theory of Bochner-Kaehler manifolds. We study the class of
Bochner-Kaehler manifolds whose scalar distribution is of special type.
Taking into account that any manifold of this class locally is a
warped product Kaehler manifold, we describe all warped product
Bochner-Kaehler metrics. We find four families of complete metrics of
this type.
\end{abstract}
\maketitle
\section{Introduction}
\vskip 2mm In \cite{GM2} we considered K\"ahler manifolds
$(M,g,J,D)$ endowed with a $J$-invariant distribution $D$ of real
codimension 2. If $D^{\perp}$ is the 2-dimensional distribution,
orthogonal to $D$, then every holomorphic section $E(m), \; m\in
M$, determines an angle $\vartheta=\angle(E(m),D^{\perp}(m))$.

A K\"ahler manifold $(M,g,J,D)$ is said to be of quasi-constant
holomorphic sectional curvatures \cite{GM2} if its holomorphic
sectional curvatures only depend on the point $m$ and the angle
$\vartheta$. The distribution $D$ of such manifolds is of pointwise
constant holomorphic sectional curvatures $a(m)$.

In \cite{GM2} we studied the case of a K\"ahler manifold of
quasi-constant holomorphic sectional curvatures when the
distribution $D$ is non-involutive. This implies $da\neq 0$ and the
1-form $\displaystyle{\eta=\frac{da}{\Vert a \Vert}}$ generates
an involutive distribution $\Delta$ \,(determined by the nullity
spaces of $\eta$). Then $D$  is a $B_0$-distribution with a geometric
function $k\neq 0$. We proved that the integral submanifolds of
$\Delta$ are $\alpha$-Sasakian manifolds of constant
$\varphi$-holomorphic sectional curvatures of type determined by
$sign \,(a+k^2)$. This result establishes a relation between
the K\"ahler manifolds of quasi-constant holomorphic sectional
curvatures and the $\alpha$-Sasakian space forms.

In this paper we show that by using $\alpha_0$-Sasakian space forms
we can construct locally all K\"ahler manifolds of quasi-constant
holomorphic sectional curvatures with $B_0$-distribution.

In Section 3 we consider warped product manifolds using as a base an
$\alpha_0$-Sasakian manifold and introduce geometrically determined
complex structure $J$. We prove that the warped product manifold with
the complex structure \,$J$ \,becomes a locally conformal K\"ahler manifold.
Further we define a K\"ahler metric on these manifolds.

In Proposition \ref{P:3.2} we prove that any warped product K\"ahler
manifold is of quasi-constant holomorphic sectional curvatures if and
only if the underlying manifold is an $\alpha_0$-Sasakian space form.

The basic result in Section 3 (Theorem \ref{T:3.1}) states as follows:
\vskip 2mm
{\it Any K\"ahler manifold of quasi-constant holomorphic sectional
curvatures with $B_0$-distribution, satisfying one of the conditions
$$a+k^2>0, \quad a+k^2=0, \quad a+k^2<0,$$
locally has the structure of a warped product K\"ahler manifold from
the Examples (E).}
\vskip 2mm
In Section 4 we study the local theory of Bochner-K\"ahler manifolds
with respect to the $J$-invariant scalar distribution generated by
the scalar curvature $\tau$ of the manifold. We prove that the
(1,0)-part of the vector field $grad\,\tau$ is a holomorphic vector
field and the function
$${\B}=\Vert \rho \Vert^2 - \frac{\tau^2}{2(n+1)}
+ \frac{\Delta \tau}{n+1}$$
is a constant, which we call the Bochner constant of the manifold.

In Section 5 we investigate Bochner-K\"ahler manifolds whose scalar
distribution is a $B_0$-distribution. In Proposition \ref{P:5.1}
we give the following curvature characterization of these manifolds:
\vskip 2mm
{\it Let $(M,g,J) \; (\dim M = 2n \geq 6)$ be a Bochner-K\"ahler manifold
with $d\tau \neq 0$. If the scalar distribution $D_{\tau}$ is a
$B_0$-distribution, then the manifold is of quasi-constant holomorphic
sectional curvatures and}
$$R = a\pi + b\Phi, \quad b \neq 0.$$
\vskip 2mm
On any of the above manifolds the function
${\bb}_0=\displaystyle{\frac{2a-b}{2}}$ is a constant. The constants
${\B}$ and ${\bb}_0$ determine uniquely the constant ${\K}$.

Theorem \ref{T:3.1} shows that the study of Bochner-K\"ahler manifolds
whose scalar distribution is a $B_0$-distribution is equivalent to
the study of warped product Bochner-K\"ahler manifolds.

In Section 6 we describe completely all families of warped product
Bochner-K\"ahler metrics in terms of the underlying $\alpha_0$-Sasakian
space form $Q_0$ and the constants ${\K},\,{\bb}_0$.

Finally we show that four of these families consist of complete metrics.

\section{Preliminaries}
\vskip 2mm Let $(M,g,J,D)$ be a $2n$-dimensional K\"ahler manifold
with metric $g$, complex structure $J$ and $J$-invariant
distribution $D$ of codimension 2. The Lee algebra of all
$\mathcal{C}^{\infty}$ vector fields on $M$ will be denoted by
${\X}M$ and $T_m M$ will stand for the tangent space to $M$ at an
arbitrary point $m \in M$. Any tangent space to $M$ has the
structure $T_m M=D(m)\oplus D^{\perp}(m),$ where $D^{\perp}(m)$ is
the 2-dimensional $J$-invariant orthogonal complement to the space
$D(m),\; m\in M$. Then the structural group of these manifolds is
the subgroup $U(n-1)\times U(1)$ of $U(n)$.

Locally we can always choose a unit vector field $\xi$ such that
$D^{\perp}=span \{\xi,J\xi\}$. The 1-forms corresponding to the
vector fields $\xi$ and $J\xi$ with respect to the metric $g$
are given by:
$$\eta(X)=g(X,\xi),\quad \tilde\eta(X)=g(X,J\xi)=-\eta(JX),
\quad X \in {\X}M.$$

The K\"ahler form\, $\Omega$ \,of the K\"ahler structure $(g,J)$
is given by\, $\Omega(X,Y)=g(JX,Y),\quad X, Y \in {\X}M$.

Let $\nabla$ be the Levi-Civita connection of the metric $g$.
The Riemannian curvature tensor $R$, the Ricci tensor $\rho$
and the scalar curvature $\tau$ of $\nabla$ are given as follows:
$$R(X,Y)Z = \nabla_X\nabla_YZ - \nabla_Y\nabla_XZ - \nabla_{[X,Y]}Z,$$
$$R(X,Y,Z,U) = g(R(X,Y)Z,U), \quad X,Y,Z,U \in {\X}M,$$
$$\rho (Y,Z) = \displaystyle{\sum_{i=1}^{2n}} R(e_i, Y, Z, e_i),
\quad Y,  Z \in T_mM,$$
$$\tau = \displaystyle{\sum_{i = 1}^{2n}\rho(e_i,e_i)},$$
where $\{e_i\}, \, i = 1,...,2n$ is an orthonormal basis for $T_m
M, \; m \in M.$ The structure $(g,J,D)$ also gives rise to the
following functions
$$\varkappa:=R(\xi,J\xi,J\xi,\xi), \quad
\sigma:=\rho(\xi,\xi)=\rho(J\xi,J\xi)$$
and tensors

$$\begin{array}{ll}
4\pi (X,Y)Z &:= g(Y,Z)X - g(X,Z)Y - 2g(JX,Y)JZ\\
[2mm]
&+ g(JY,Z)JX - g(JX,Z)JY,\\
[2mm]
8\Phi(X,Y)Z&:=g(Y,Z)(\eta(X)\xi+\tilde\eta(X)J\xi)
-g(X,Z)(\eta(Y)\xi+\tilde\eta(Y)J\xi)\\
[2mm]
&+g(JY,Z)(\eta(X)J\xi-\tilde\eta(X)\xi)-g(JX,Z)(\eta(Y)J\xi
-\tilde\eta(Y)\xi)\\
[2mm]
&-2g(JX,Y)(\eta(Z)J\xi-\tilde\eta(Z)\xi)\\
[2mm]
&+(\eta(Y)\eta(Z)+
\tilde\eta(Y)\tilde\eta(Z))X
-(\eta(X)\eta(Z)+\tilde\eta(X)\tilde\eta(Z))Y\\
[2mm]
&+(\eta(Y)\tilde\eta(Z)-\tilde\eta(Y)\eta(Z))JX
-(\eta(X)\tilde\eta(Z)-\tilde\eta(X)\eta(Z))JY\\
[2mm]
&-2(\eta(X)\tilde\eta(Y)-\tilde\eta(X)\eta(Y))JZ,\\
[2mm]
\Psi (X,Y)Z &:= \eta(Y)\eta(Z)\tilde\eta(X)J\xi
- \eta(X)\eta(Z)\tilde\eta(Y)J\xi\\
[2mm]
&+ \eta(X)\tilde\eta(Y)\tilde\eta(Z)\xi -
\eta(Y)\tilde\eta(X)\tilde\eta(Z)\xi\\
[2mm]
&=(\eta \wedge\tilde\eta)(X, Y)(\tilde\eta(Z)\xi
- \eta(Z)J\xi),\quad X,Y,Z \in {\X}M.
\end{array}\leqno(2.1)$$
These basic tensors are invariant under the action of the
structural group $U(n-1)\times U(1)$ in the sense of \cite{TV}.

A K\"ahler manifold  $(M,g,J,D) \; (\dim M=2n\geq 4)$ with
$J$-invariant distribution $D$ is of quasi-constant holomorphic
sectional curvatures (briefly a QCH-manifold) if and only if
\cite{GM2}
$$R=a\pi+b\Phi+c\Psi, \leqno(2.2)$$
where $\pi, \,\Phi,\, \Psi$ are the tensors (2.1) and
$$ \begin{array}{l}
\displaystyle{ a = \frac{\tau - 4 \sigma + 2\varkappa}{n(n-1)},
\quad b = \frac{- 2\tau +4(n+2)\sigma - 4(n+1)\varkappa}{n(n-1)}},\\
[4mm]
\displaystyle{c = \frac{\tau - 4(n+1)\sigma +
(n+1)(n+2)\varkappa}{n(n-1)}}.
\end{array}\leqno(2.3)$$

Let $(M,g)$ be a Riemannian manifold endowed with a unit vector
field $\xi$, and $\eta$ be the 1-form corresponding to $\xi$
with respect to the metric $g$. We assume that the distribution
$\Delta$, determined by the nullity spaces of $\eta$, is involutive.
An involutive distribution $\Delta$ is characterized by the
condition \cite{GM1}
$$d\eta=-\theta \wedge \eta,$$
where $\theta(X)=g(\nabla_{\xi}\,\xi,X),\quad X \in {\X}M.$

A $\mathcal{C}^{\infty}$ function $u$ on $M$ is said to be a
proper function for the involutive distribution $\Delta$ if
\cite{GM1}
$$du=\xi(u)\,\eta=\Vert du \Vert \,\eta.$$

Now, let $(M,g,J)$ be a K\"ahler manifold endowed with a unit
vector field $\xi$ and $\eta, \, \tilde\eta$ be the 1-forms
corresponding to the vector fields $\xi, \, J\xi$, respectively.
Then the vector field $\xi$ generates the distributions:
$$\Delta(m):=\{x \in T_m M\, \vert\, \eta(x)=0\},\quad m \in M,$$
$$D(m)=\{x_0 \in T_m M\; \vert\; \eta(x_0)=\tilde\eta(x_0)=0\},\quad
D^{\perp}(m)= span \{\xi, J\xi\}, \quad m \in M.$$

As a rule we shall use the following denotations for the
different kinds of vectors (vector fields) on $M$:
$$X \in T_m M \;({\X}M),\quad x \in \Delta \; ({\X}\Delta),
\quad x_0 \in D(m) \;({\X}D), \quad m \in M.$$

Any vector $X \in T_m M$ can be decomposed in a unique way in the
form
$$X=x_0+\tilde\eta(X)J\xi+\eta(X)\xi,$$
where $x_0$ is the projection of $X$ into $D(m)$.

The following notion is essential in our considerations:
\begin{defn}\label{D:2.1} \cite{GM2, GM3}
A $J$-invariant distribution $D$ generated by the unit vector
field $\xi$ is said to be a $B_0$-distribution if it satisfies
the following conditions:
$$\begin{array}{l}
i) \displaystyle{ \quad \nabla _{x_0}\,\xi =
\frac{k}{2}\,x_0, \quad k\neq 0,}\quad x_0 \in D,\\
[2mm]
ii) \quad \nabla_{J\xi}\,\xi=-p^*J\xi,\\
[2mm]
iii)\quad \nabla _{\xi}\,\xi = 0. \\
[2mm]
\end{array}\leqno(2.4)$$
\end{defn}

The conditions (2.4) in Definition 2.1 are equivalent to
the equality
$$\nabla_X\,\xi=\frac{k}{2}\{X-\eta(X)\xi-\tilde\eta(X)J\xi\}
-p^*\tilde\eta(X)J\xi, \quad X\in {\X}M. \leqno (2.5)$$

The last formula implies that
$$dk = \xi(k)\,\eta, \quad p^* = -\frac{\xi(k)+k^2}{k}.\leqno{(2.6)}$$

Introducing the relative divergences $div_0\,\xi$ and
$div_0\,J\xi$ (relative codifferentials $\delta_0\eta$ and
$\delta_0\tilde\eta$) of the vector fields $\xi$ and $J\xi$
(1-forms $\eta$ and $\tilde\eta$) with respect to the distribution
$D$ by
$$div_0\,\xi = -\delta_0\eta =
\sum_{i = 1}^{2(n-1)}(\nabla_{e_i}\eta)e_i, \quad div_0\,J\xi =
-\delta_0\tilde\eta = \sum_{i =
1}^{2(n-1)}(\nabla_{e_i}\tilde\eta)e_i,$$ $\{e_1,...,e_{2(n-1)}\}$
being an orthonormal basis of $D(m), \, m \in M$, from (2.5) we
obtain:
$$k=\frac{div_0\,\xi}{n-1},\quad div_0\,J\xi=0.$$

Let $(M,g,J,D) \; (\dim M=2n\geq 6)$ be a K\"ahler QCH-manifold
with non-involutive distribution $D
\;(D^{\perp}=span\{\xi,J\xi\}).$ From (2.5) it follows that $D$ is
a $B_0$-distribution with geometric function $k \neq 0$ and the
equality (2.2) implies that $D$ is of pointwise constant
holomorphic sectional curvatures $a$. In \cite {GM2} we have shown
that the function $a+k^2$ plays an important role in the study of
the class of the K\"ahler QCH-manifolds with $B_0$-distribution.
We divided this class into three subclasses with respect to $sign
\, (a+k^2)$ \cite {GM3}:
$$a+k^2>0, \quad a+k^2=0, \quad a+k^2<0.$$

In \cite{GM3} we have proved that the integral submanifolds of the
distribution $\Delta$ are $\alpha$-Sasakian space forms.

Next we give some basic notions concerning the class of
$\alpha$-Sasakian manifolds.

Let $(Q_0,g_0,\varphi,\tilde \xi_0,\tilde\eta_0)\,(\dim Q_0 =2n-1 \geq 5)$
be an almost contact Riemannian manifold. The structures of the manifold
satisfy the following conditions:
$$\begin{array}{l}
\tilde \eta_0(x)=g_0(\tilde \xi_0,x),\quad \varphi\tilde\xi_0=0,\quad
{\varphi}^2 x=-x+\tilde\eta_0(x)\tilde\xi_0,\quad
x \in {\X}Q_0,\\
[2mm]
g_0(\varphi x, \varphi y)=g_0(x, y)-\tilde\eta_0(x)\tilde\eta_0(y),
\quad x, y \in {\X}Q_0.
\end{array} \leqno (2.7)$$
Let $\mathcal{D}^0$ be the Levi-Civita connection of the metric $g_0$.
The manifold $(Q_0,g_0,\varphi,\tilde \xi_0,\tilde\eta_0)$ is
called an $\alpha_0$-Sasakian manifold {\cite {JV}} if
$$ \mathcal{D}^0_x\,\tilde \xi_0 = \alpha_0\;\varphi x, \quad
x \in {\X}Q^{2n-1}_0, \quad \alpha_0 = const;\leqno (2.8)$$
$$(\mathcal{D}^0_x\varphi)y =
\alpha_0\{\tilde \eta_0(y)x-g_0(x,y)\tilde\xi_0\},
\quad x, y \in {\X}Q_0.\leqno (2.9)$$
From (2.8) it follows that
$$d\tilde\eta_0(x,y)=2\alpha_0\,g_0(\varphi x,y),\quad
x, y \in {\X}Q_0.$$
If $K^0$ is the Riemannian curvature tensor of $g_0$, then
(2.7), (2.8) and (2.9) imply that
$$K^0(x, y)\tilde \xi_0=\alpha^2_0\{\tilde\eta_0(y)x-\tilde\eta_0(x)y\},
\quad x, y \in {\X}Q_0.\leqno (2.10)$$
The last equality gives the geometric meaning of the constant $\alpha_0$:
$$K^0(x_0,\tilde\xi_0,\tilde\xi_0,x_0)=\alpha^2_0,\quad
x_0 \perp \tilde\xi_0, \quad g_0(x_0,x_0)=1.$$
The change of the direction of $\tilde\xi_0$ has as a consequence
the change of $sign \,\alpha_0$. Thus we can assume that
$\alpha_0>0.$

The class of the Sasakian manifolds is the subclass of
$\alpha_0$-Sasakian manifolds with $\alpha_0=1$.

The distribution orthogonal to $\tilde\xi_0$ is denoted by $D$:
$$D(m):=\{x_0 \in T_m Q_0 \,\vert\, \tilde\eta_0(x_0)=0\},\quad
m \in Q_0.$$ Then the vector space $(D,\varphi)$ is a Hermitian
vector space at every point of $Q_0$.

An $\alpha_0$-Sasakian manifold is said to be an
$\alpha_0$-Sasakian space form {\cite {O}} if it is of constant
$\varphi$-holomorphic sectional curvatures $H_0$, i.e.
$$K^0(x_0,\varphi x_0,\varphi x_0,x_0)=H_0,\quad x_0 \in D,
\quad g_0(x_0,x_0)=1.$$

Any $\alpha_0$-Sasakian space form is characterized by the following
curvature identity {\cite {O,JV}}
$$K^0=\displaystyle{\frac{H_0+3{\alpha}^2_0}{4}}\, \pi_{01}
+\displaystyle{\frac{H_0-{\alpha}^2_0}{4}}\,(\pi_{02}-\pi_{03}),
\leqno(2.11)$$
where $\pi_{01}, \pi_{02}, \pi_{03}$ are the basic invariant tensors
with respect to the structure $(g_0,\tilde\xi_0,\tilde\eta_0)$:
$$\begin{array}{ll}
\pi_{01}(x,y,z,u)=&g_0(y,z)g_0(x,u)-g_0(x,z)g_0(y,u),\\
[2mm]
\pi_{02}(x,y,z,u)=&g_0(\varphi y,z)g_0(\varphi x,u)-
g_0(\varphi x,z)g_0(\varphi y,u)-2g_0(\varphi x,y)g_0(\varphi z,u),\\
[2mm]
\pi_{03}(x,y,z,u)=&g_0(y,z)\tilde\eta_0(x)\tilde\eta_0(u)-
g_0(x,z)\tilde\eta_0(y)\tilde\eta_0(u)\\
[2mm]
&+\tilde\eta_0(y)\tilde\eta_0(z)g_0(x,u)-
\tilde\eta_0(x)\tilde\eta_0(z)g_0(y,u),
\end{array} \leqno (2.12) $$
for any $x,y,z,u \in {\X}Q_0$.

Keeping in mind the projections $g_0-\tilde\eta_0\otimes\tilde\eta_0$
and $\tilde\eta_0\otimes\tilde\eta_0$ of the metric $g_0$ onto the
distribution $D$ and $span\{\tilde\xi_0\}$, we consider the following
{\it bihomothetical} transformation of the structure
$(g_0, \tilde\xi_0, \tilde\eta_0)$:
$$\begin{array}{ll}
g&=p^2\{(g_0-\tilde\eta_0\otimes\tilde\eta_0)+
q^2\,\tilde\eta_0\otimes\tilde\eta_0\}\\
[2mm]
&=p^2\{g_0+(q^2-1)\tilde\eta_0\otimes \tilde\eta_0\},\\
[2mm]
\tilde\xi&=\displaystyle{\frac{1}{pq}}\,\tilde\xi_0, \quad
\tilde\eta =pq\,\tilde\eta_0,
\end{array} \leqno(2.13)$$
where $p$ and $q$ are positive constants.
In the case $p=q$, the transformation (2.13) is called a
D-homothetic deformation of the structure \cite {T}.

If $\mathcal{D}$ is the Levi-Civita connection  of the new metric $g$,
then
$$\begin{array}{l}
\mathcal{D}_x\,\tilde\xi=\displaystyle{\alpha_0\frac{q}{p}}\,\varphi x,\\
[2mm]
(\mathcal{D}_x\varphi)y=\displaystyle{\alpha_0\frac{q}{p}}\,
\{\tilde\eta(y)x-g(x,y)\tilde\xi\},
\end{array}\leqno(2.14)$$
where $x, y \in {\X}Q_0$.

From (2.14) it follows that $(Q_0,g,\varphi,\tilde\xi,\tilde\eta)$
is an $\alpha$-Sasakian manifold with
$\alpha=\displaystyle{\frac{q}{p}}\,\alpha_0.$

If $\pi_1, \pi_2, \pi_3$ denote the basic invariant tensors (2.12) with
respect to the structure $(g,\tilde\xi,\tilde\eta)$ given by (2.13),
then the relations between $\pi_1, \pi_2, \pi_3$ and
$\pi_{01}, \pi_{02}, \pi_{03}$ are the following:
$$\begin{array}{l}
\pi_1=p^4\{\pi_{01}+(q^2-1)\pi_{03}\},\\
[2mm]
\pi_2=p^4\,\pi_{02},\\
[2mm]
\pi_3=p^4q^2\,\pi_{03}.
\end{array}\leqno(2.15)$$

The Levi-Civita connection $\mathcal{D}$ of $g$, as a consequence of
(2.13), is expressed by the equality
$$\mathcal{D}_x \,y=\mathcal{D}^0_x\,y+
\alpha_0(q^2-1)\{\tilde\eta_0(y)\, \varphi x+
\tilde\eta_0(x)\, \varphi y\},\quad x, y \in {\X}Q_0.$$
Then the curvature tensor $K$ of the metric $g$ has the form:
$$K=\frac{H+3\alpha^2}{4}\,\pi_1+
\frac{H-\alpha^2}{4}\,(\pi_2-\pi_3),\leqno (2.16)$$
where
$$H=\frac{H_0+3\alpha^2_0}{p^2}-
\frac{3q^2}{p^2}\,\alpha^2_0, \quad \alpha=\frac{q}{p}\,\alpha_0.
\leqno (2.17)$$
Hence
$$H+3\alpha^2=\frac{H_0+3\alpha^2_0}{p^2},
\leqno (2.18)$$
which shows that $sign \,(H_0+3\alpha^2_0)$ is invariant under
the bihomothetical changes (2.13) of the structure
$(g_0,\tilde\xi_0,\tilde\eta_0)$.

There are three types of $\alpha_0$-Sasakian space forms
corresponding to the $sign$ of the constant $H_0+3\alpha^2_0$:
$$H_0+3\alpha^2_0>0,\quad H_0+3\alpha^2_0=0,\quad H_0+3\alpha^2_0<0.$$
The formula (2.18) shows that the bihomothetical transformations
(2.13) of the structure $(g_0,\tilde\xi_0,\tilde\eta_0)$ do
not change the type of the space form.

\section{Warped Product K\"ahler Manifolds}
\vskip 2mm
Now let $(Q_0, g_0,\varphi,\tilde\xi_0,\tilde\eta_0)\;
(\dim Q_0= 2n-1\geq 5)$ be an
$\alpha_0$-Sasakian manifold with constant $\alpha_0>0$. Further,
let $\mathbb{R}$ be equipped with a coordinate system $O\textbf{e}$ and
the standard inner product determined by $\textbf{e}^2=1$. For an
arbitrary point $t \in \mathbb{R}$ we denote
$\displaystyle{\bar\xi:=\frac{d}{d{t}}}$ and $\bar\eta:=d{t}$.

The manifolds $Q_0$ and $\mathbb{R}$ generate the product
manifold $N=Q_0\times \mathbb{R}$ with the standard product
metric $g_0+ \bar\eta\otimes\bar\eta$.

Let $p(t),\, t\in I \subset \mathbb{R}$ with $p(t)>0, \,p(t_0)=1$
be a real function. On any leaf
$Q(t)\; (\dim Q=2n-1\geq 5)$ \,($t$ - fixed)
with the $\alpha_0$-Sasakian structure
$(g_0,\varphi,\tilde\xi_0,\tilde\eta_0)$, we define the structure
$(\bar g,\varphi,\tilde{\bar\xi},\tilde{\bar\eta})$:
$$\bar g(t):=p^2(t)g_0, \quad \tilde{\bar\xi}:=
\displaystyle{\frac{1}{p(t)}\,}\tilde\xi_0,
\quad \tilde{\bar\eta}:=p(t)\,\tilde\eta_0.\leqno(3.1)$$
This structure is homothetical to $(g_0,\varphi,\tilde\xi_0,
\tilde\eta_0)$ with constant $p(t)$ ($t$ - fixed).
Then the structure
$(\bar g,\varphi,\tilde{\bar\xi},\tilde{\bar\eta})$ defined by (3.1)
\,($t$-fixed) is an $\bar\alpha$-Sasakian with
$$\bar\alpha(t)=\displaystyle{\frac{1}{p(t)}\,\alpha_0}.$$

We furnish the manifold $N=Q_0 \times \mathbb{R}$ with the
warped product metric \cite{BO}:
$$G:=p^2(t)g_0+\bar\eta\otimes\bar\eta.\leqno (3.2)$$

Then $\bar\eta$ and $\tilde{\bar\eta}$ are the 1-forms corresponding
to $\bar\xi$ and $\tilde{\bar\xi}$ with respect to the metric $G$:
$$\bar\eta(X)=G(\bar\xi,X),\quad \tilde{\bar\eta}(X)=
G(\tilde{\bar\xi},X),\quad X \in {\X}N.$$ Let $T_m N$ be the
tangent space to $N$ at an arbitrary point $m \in N$. The
structure $(G,\bar\xi,\tilde{\bar\xi})$ gives rise to the
following distributions on $N$:
$$\begin{array}{l}
D(m):=\{x_0 \in T_m N \,\vert\, \tilde{\bar\eta}(x_0)=
\bar\eta(x_0)=0\},\quad D^{\perp}(m):=
span \,\{\bar\xi,\tilde{\bar\xi}\},\\
[2mm] \Delta(m):=\{x \in T_m N \,\vert\, \bar\eta(x)=0\},\quad m
\in N.
\end{array}$$
Then any vector (vector field) $X$ in $T_m N \, \,({\X}N)$ can be
decomposed uniquely as follows:
$$X=x_0+\tilde{\bar\eta}(X)\tilde{\bar\xi}+\bar\eta(X)\bar\xi,$$
where $x_0$ is the projection of $X$ into $D$.

Any vector (vector field) $x$ in $\Delta \, \,({\X} \Delta)$ can be
decomposed in a unique way as follows:
$$x=x_0+\tilde{\bar\eta}(x)\tilde{\bar\xi},$$
where $x_0$ is the projection of $x$ into $D$.

Denoting by $\bar\nabla$ and $\bar{\mathcal{D}}$ the Levi-Civita
connections of the manifold $(N,G)$ and the submanifold
$(Q(t),\bar g(t))$,
respectively, we have
\begin{lem}\label{L:3.1}
The Gauss and Weingarten formulas for the submanifold
$(Q(t),\bar g(t))$
of $N$ with normal vector field $\bar\xi$ are:
$$\bar\nabla_x\,y=\bar{\mathcal{D}}_x\,y -\bar\xi(\ln{p})\,\bar g(x,y)\,
\bar\xi, \quad x, y \in {\X}\Delta,\leqno(3.3)$$
$$\bar\nabla_x\,\bar\xi=\bar\xi(\ln{p})\,x,
\quad x \in {\X}\Delta.\leqno(3.4)$$
\end{lem}
{\it Proof}. Since $N=Q_0\times\mathbb{R}$, then \cite{BO}
$$\bar\nabla_{\bar\xi}\,\bar\xi=0 $$
and $[\bar\xi,\tilde\xi_0]=[\bar\xi,x_0]=0$ for every $x_0 \in {\X}D.$
From (3.1) it follows
$$[\bar\xi,\tilde{\bar\xi}]=\bar\nabla_{\bar\xi}\,\tilde{\bar\xi}-
\bar\nabla_{\tilde{\bar\xi}}\,\bar\xi=
\displaystyle{\frac{1}{p}}\,[\bar\xi,\tilde\xi_0]
-\bar\xi(\ln{p})\,\tilde{\bar\xi}=
-\bar\xi(\ln{p})\,\tilde{\bar\xi}.$$

Using the Koszul formula for the Levi-Civita connection $\bar\nabla$
of the metric $G$
$$\begin{array}{ll}
2G(\bar\nabla_XY,Z)&=XG(Y,Z)+YG(Z,X)-ZG(X,Y)\\
[2mm]
&+ G([X,Y],Z)-G([Y,Z],X)+G([Z,X],Y),\quad X, Y, Z \in {\X}N,
\end{array}$$
we calculate the components of $\bar\nabla_XY,\quad X, Y \in {\X}N$:
$$\begin{array}{l}
\bar\nabla_{x_0}\,\tilde{\bar\xi}_{\vert \Delta}
=\bar{\mathcal{D}}_{x_0}\,\tilde{\bar\xi} =\bar\alpha\,\varphi x_0,\quad
\bar\nabla_{\tilde{\bar\xi}}\,{x_0}_{\vert \Delta}
=\bar{\mathcal{D}}_{\tilde{\bar\xi}}\,x_0 =\bar\alpha\,\varphi x_0,\quad
x_0 \in {\X}D,\\
[2mm]
\bar\nabla_{\tilde{\bar\xi}}\,\tilde{\bar\xi} _{\vert \Delta}=0,\quad
\bar\nabla_{\tilde{\bar\xi}}\,\tilde{\bar\xi} =-\bar\xi(\ln {p})\,\bar\xi,
\quad \bar\nabla_{\tilde{\bar\xi}}\,\bar\xi=\bar\xi(\ln {p})
\,\tilde{\bar\xi},\quad \bar\nabla_{\bar\xi}\,\tilde{\bar\xi}=0,\\
[2mm]
\bar\nabla_{x_0}\,\bar\xi=\bar\nabla_{\bar\xi}\, x_0=\bar\xi(\ln p)\,x_0,
\quad G(\bar\nabla_{x_0}\,y_0,\bar\xi)=-\bar\xi(\ln p)\,\bar g(x_0,y_0),
\quad x_0, y_0 \in {\X}D.
\end{array} \leqno(3.5)$$
The above equalities imply the assertion of the lemma.
\hfill{\bf QED}

\vspace{2mm}
Next we define a complex structure $J$ on the
warped product manifold $(N,G)$ by the help of the geometric
structures $(\varphi,\tilde{\bar\xi},\bar\xi)$ in the following
way:
$$J_{\vert D}:=\varphi,\quad J\bar\xi:=\tilde{\bar\xi},\quad
J\tilde{\bar\xi}:=-\bar\xi.\leqno(3.6)$$
Thus $(N,G,J)$ becomes an almost Hermitian manifold.

More precisely we have
\begin{lem} \label{L:3.2}
The covariant derivative of the structure $J$ on $(N,G)$, given
by $(3.6)$, satisfies the following identity:
$$(\bar\nabla_X J)Y=\{\bar\xi(\ln p)-\bar\alpha\}\{G(X,Y)\tilde{\bar\xi}
+G(JX,Y)\bar\xi-\tilde{\bar\eta}(Y)X-\bar\eta(Y)JX\},\quad
X,Y \in {\X}N.\leqno(3.7)$$
\end{lem}

{\it Proof.} Let $X,Y \in {\X}N,\; x,y$ be the projections of
$X,Y$ into ${\X}\Delta$ and $x_0, y_0$ - the projections of
$X,Y$ into ${\X}D$. Then we have
$$\begin{array}{l}
\vspace{2mm}
X=x+\bar\eta(X)\bar\xi=x_0+\tilde{\bar\eta}(X)\tilde{\bar\xi}+
\bar\eta(X)\bar\xi,\\
\vspace{2mm}
Y=y+\bar\eta(Y)\bar\xi=y_0+\tilde{\bar\eta}(Y)\tilde{\bar\xi}+
\bar\eta(Y)\bar\xi;\end{array}
\leqno{(3.8)}$$
$$Jy=\varphi y-\tilde{\bar\eta}(y)\bar\xi. \leqno(3.9)$$
Using (3.3), (3.6) and (3.9) we calculate
$\bar\nabla_x\,Jy-J\bar\nabla_x\,y:$
$$(\bar\nabla_x J)y=\{\bar\xi(\ln p)-\bar\alpha\}
\{\bar g(x,y)\tilde{\bar\xi}+\bar g(\varphi x,y)\bar\xi
-\tilde{\bar\eta}(y)x\},\quad x,y \in {\X}\Delta.
\leqno (3.10)$$
Using (3.3), (3.4), (3.6) and (3.8) we compute
$\bar\nabla_x\,J\bar\xi-J\bar\nabla_x\,\bar\xi$:
$$(\bar\nabla_x J)\bar\xi=-\{\bar\xi(\ln p)-\bar\alpha\}\,\varphi x,
\quad x \in {\X}\Delta. \leqno(3.11)$$

Finally, by using (3.5) and (3.6), we find
$$(\bar\nabla_{\bar\xi} J)\tilde{\bar\xi}=0,\quad
(\bar\nabla_{\bar\xi} J)\bar\xi=0,
\quad (\bar\nabla_{\bar\xi}\,J)x_0=0, \quad x_0 \in {\X}D,$$
which imply
$$(\bar\nabla_{\bar\xi} J)X=0, \quad X \in {\X}N.\leqno (3.12)$$
Now (3.10), (3.11) and (3.12) give (3.7).
\hfill{\bf QED}

\vspace {2mm}
The identity (3.7) implies that the almost complex structure $J$,
defined by (3.6) is integrable, i.e. $(N,G,J)$ is a Hermitian manifold.
More precisely, this manifold is in the class $W_4$, according to
the classification of almost Hermitian manifolds given in \cite{GH},
with Lee form $\{\bar\alpha-\bar\xi(\ln p)\}\,\bar\eta.$
Since $\bar\alpha-\bar\xi(\ln p)$ is a function of $t$
and $\bar\eta = d{t}$, then
the Lee form is closed in all dimensions $2n\geq4$ and $(N,G,J)$ is
a locally conformal K\"ahler manifold.

Let $q(t)>0,\, t\in I \subset {\R}$ be a real
$\mathcal{C}^{\infty}$-function satisfying the condition $q(t_0)=1$.
We consider the following complex dilatational transformation of the
structure $(G,\bar\xi)$ (cf \cite{GM2}):
$$\begin{array}{l}
\vspace{2mm}
\textsl{g}=G+(q^2(t)-1)(\bar\eta\otimes\bar\eta+
\tilde{\bar\eta}\otimes\tilde{\bar\eta}),\\
\vspace{2mm}
\xi=\displaystyle{\frac{1}{q}\,\bar\xi}, \quad \eta=q\,\bar\eta;
\quad \tilde\xi=\displaystyle{\frac{1}{q}\,\tilde{\bar\xi}},
\quad \tilde\eta=q\,\tilde{\bar\eta}.
\end{array}\leqno{(3.13)}$$

\begin{prop} \label{P:3.1}
The Hermitian manifold $(N,\textsl{g},J,\xi, \tilde\xi)$ with
structures given by $(3.6)$ and $(3.13)$ is
K\"ahlerian if and only if the positive functions
$p(t)$ and $q(t)$ are related as follows
$$q^2(t)=\frac{1}{\alpha_0}\,\bar\xi(p).$$
\end{prop}

{\it Proof.} From (3.13) we obtain that the K\"ahler form
$\textsl{g}(JX,Y), \,X,Y \in {\X}N$, is closed if and only if
$\bar\xi(\ln p)-\bar\alpha q^2=0$.  Since
$\bar\alpha(t)=\displaystyle{\frac{1}{p(t)}\,\alpha_0}$,
then $\displaystyle{q^2=\frac{\bar\xi(p)}{\alpha_0}}$.
\hfill{\bf QED}
\vskip 2mm
\begin{rem}\label{R:3.1}
In view of (3.13) the 1-form $\eta$ is closed, i.e. $\eta=ds$ locally.
Taking into account the equality $\eta=q\,\bar\eta$ and choosing
$s(t_0)=0$ we obtain
$s(t)=\displaystyle{\int_{t_0}^{t}q(t)dt}.$
Considering the functions $p(t)$ and $q(t)$ as functions of $s$
we have
$$\eta=ds, \quad  \xi =\frac{d}{ds}, \quad
q=\sqrt{\frac{1}{\alpha_0}\frac{dp}{dt}}=\frac{1}{\alpha_0}\,\xi(p)
=\frac{1}{\alpha_0}\frac{dp}{ds}=\frac{p'}{\alpha_0}.$$
\end{rem}
\vskip 2mm

From now on the derivatives with respect to the parameter $s$
will be denoted as usual by $()', ()'',...$

Let $Q_0(g_0,\varphi,\tilde\xi_0,\tilde\eta_0)$ be an
$\alpha_0$-Sasakian manifold and $p(s),\, s \in I_0 \subset \mathbb{R}$,
be a $\mathcal{C}^{\infty}$-real function, satisfying the
conditions
$$p(s)>0,\quad p'(s)>0,\quad p(0)=1, \quad p'(0)=\alpha_0.\leqno{(3.14)}$$

We consider the following K\"ahler manifolds
$(N,\textsl{g},J,\xi)\, (\dim N=2n \geq 6)$:
$$\begin{array}{ll}
N=Q_0\times\mathbb{R},\\
[2mm]
\textsl{g}=\bar g- \tilde{\bar\eta}\otimes\tilde{\bar\eta}+
\displaystyle{\frac{{p'}^2}{\alpha_0^2}}\,(\bar\eta\otimes\bar\eta+
\tilde{\bar\eta}\otimes\tilde{\bar\eta})\\
[2mm]
=p^2(s)\displaystyle{\left\{g_0+
\left(\frac{{p'}^2(s)}{\alpha^2_0}-1\right)\,
\tilde\eta_0\otimes\tilde\eta_0\right\}+
\frac{{p'}^2(s)}{\alpha_0^2}\,\bar\eta\otimes\bar\eta},\\
[2mm]
\tilde{\bar\xi}=\displaystyle{\frac{1}{p}\,\tilde\xi_0,\quad
\tilde{\bar\eta}= p\,\tilde\eta_0,}\quad
\tilde\xi=\displaystyle{\frac{\alpha_0}{p'}\,\tilde{\bar\xi}=
\frac{\alpha_0}{pp'}\,\tilde\xi_0,\quad
\tilde\eta=\frac{p'}{\alpha_0}\,\tilde{\bar\eta}=
\frac{pp'}{\alpha_0}\,\tilde\eta_0},\\
[2mm]
\displaystyle{\bar\xi=\frac{d}{d{t}}},\quad \bar\eta=d{t},
\quad \xi=\displaystyle{\frac{\alpha_0}{p'}\,\bar\xi,\quad
\eta=\frac{p'}{\alpha_0}\,\bar\eta,}\\
[2mm]
J_{\vert D}:=\varphi,\quad J\xi:=\tilde\xi,\quad
J\tilde\xi:=-\xi,
\end{array}\leqno(3.15)$$
where the function $p(s),\, s \in I_0 \subset \mathbb{R},$
satisfies the conditions (3.14) and $\alpha_0=const>0$.

In what follows we call the manifolds $(N,\textsl{g},J,\xi)$
given by (3.15) {\it warped product K\"ahler manifolds}.

If $\nabla$ is the Levi-Civita connection of the K\"ahler metric
$\textsl{g}$, then the equalities (3.13), (3.4) and
$\nabla_{\xi}\,\xi=0$ imply that
$$\nabla_X\,\xi=\displaystyle{\frac{p'}{p}\{X-\tilde\eta(X)\tilde\xi-
\eta(X)\xi\}+\frac{{p'}^2+pp''}{pp'}\,\tilde\eta(X)\tilde\xi,
\quad X \in {\X}N.}\leqno(3.16)$$
Sice $\xi=\displaystyle{\frac{\alpha_0}{p'}\,\bar\xi}$, then the unit
vector field $\xi$ generates the same distributions
$\Delta,\; D$ and $D^{\perp}$. Comparing the equality (3.16) with
(2.5) we conclude that:
\vskip 1mm
{\it The distribution $D$ of the K\"ahler manifold $(3.15)$ is a
$B_0$-distribution with functions}
$$\displaystyle{\frac{k}{2}=\frac{p'}{p}, \quad
p^*=-\frac{{p'}^2+pp''}{pp'}.} \leqno(3.17)$$
Let $R,\, \rho, \, \tau$ be the curvature tensor,
the Ricci tensor and the scalar curvature of $\textsl{g}$, respectively.
From the equality (3.16) we compute:
$$\begin{array}{ll}
R(X,Y)\xi = & \displaystyle{\frac{p''}{p}\,\{\eta(X)Y-\eta(Y)X}\\
[2mm]
& -\tilde\eta(X)JY+\tilde\eta(Y)JX+2\textsl{g}(JX,Y)J\xi\}\\
[2mm] & \displaystyle{+\frac{pp'''-p'p''}{pp'}\,
\{\eta(X)\tilde\eta(Y)-\eta(Y)\tilde\eta(X)\}J\xi},\quad X,Y \in {\X}N.
\end{array}\leqno{(3.18)}$$
Taking a trace in (3.18) we find
$$\rho(X,\xi)=-\left\{\frac{pp'''-p'p''}{pp'}+2(n+1)
\frac{p''}{p}\right\}\eta(X),
\quad X \in {\X}N.\leqno{(3.19)}$$
Then the functions $\varkappa$ and $\sigma$ are
$$\varkappa =R(\xi,J\xi,J\xi,\xi)=
-\left\{\frac{pp'''-p'p''}{pp'}+4\frac{p''}{p}\right\},
\leqno{(3.20)}$$
$$\sigma=\rho(\xi,\xi)=\rho(J\xi,J\xi)=
-\left\{\frac{pp'''-p'p''}{pp'}+2(n+1)
\frac{p''}{p}\right\}.\leqno{(3.21)}$$

\begin{prop}\label{P:3.2}
Any warped product K\"ahler manifold $(N,\textsl{g},J,\xi)$ given by
$(3.15)$ is of quasi-constant holomorphic sectional curvatures if
and only if the base manifold $Q_0$ is an $\alpha_0$-Sasakian space form.
\end{prop}
{\it Proof.} Let $(N,\textsl{g},J,\xi)$ be a warped product K\"ahler
manifold defined by (3.15). Then any leaf $Q(s)$ with the restriction
of the metric $\textsl{g}$ is a submanifold of $N$ with unit
normal vector field $\xi$. Denoting by $\mathcal{D}$ the
Levi-Civita connection of the restriction of $\textsl{g}$
onto $Q(s)$ we find the Gauss formula:
$$\nabla_x\,y=\mathcal{D}_x\,y -
\left(\displaystyle{\frac{p'}{p}\,\textsl{g}(x,y)+
\frac{p''}{p'}\,\tilde\eta(x)\tilde\eta(y)}\right)\xi,\quad
x, y \in {\X}\Delta.\leqno(3.22)$$

The submanifold $Q(s)$ carries the standard induced structure
$(\textsl{g},\varphi,\tilde\xi,\tilde\eta)$ \cite{T1,T2}, where
$$\varphi x=Jx+\tilde\eta(x)\xi, \quad x\in {\X}\Delta.\leqno{(3.23)}$$
Then the Gauss equation is
$$R_{\vert \Delta}= K-\displaystyle{\frac{{p'}^2}{p^2}\,
{\pi_1}_{\vert \Delta}-\frac{p''}{p}\,{\pi_3}_{\vert \Delta}},
\leqno(3.24)$$
where $K$ is the curvature tensor of $Q(s)$ and $\pi_1, \pi_3$
are the tensors (2.12) with respect to the structure
$(\textsl{g},\varphi,\tilde\xi,\tilde\eta)$ on $Q(s)$.

The equalities (3.22) and (3.16) imply that
$$\begin{array}{l}
\mathcal{D}_x\,\tilde\xi=\displaystyle{\frac{p'}{p}}\,\varphi x, \quad
x \in {\X}\Delta,\\
[2mm]
(\mathcal{D}_x\varphi)y=\displaystyle{\frac{p'}{p}}\,
\{\tilde\eta(y)x-\textsl{g}(x,y)\tilde\xi\}, \quad x,y \in {\X}\Delta.
\end{array}$$
Hence $Q(s)$ with the induced almost contact Riemannian structure
is an $\alpha$-Sasakian manifold with
$$\alpha(s)=\displaystyle{\frac{p'(s)}{p(s)}=const>0,
\quad s - \rm{fixed}}.\leqno(3.25)$$

Now let the base manifold $(Q_0,g_0,\varphi,\tilde \xi_0,\tilde\eta_0)$
be an $\alpha_0$-Sasakian space form with constant $\varphi$-holomorphic
sectional curvatures $H_0$.

The equalities (3.15) imply that the restriction of the
K\"ahler metric $\textsl{g}$ onto the leaf $Q(s),
\, s \in I_0 \subset \mathbb{R}$ \; ($s$-fixed), is
$$\textsl{g}(s)=p^2(s)\displaystyle{\left\{g_0+
\left(\frac{{p'}^2(s)}{\alpha^2_0}-1\right)
\, \tilde\eta_0\otimes\tilde\eta_0\right\}},$$
which means that the metrics $\textsl{g}(s)$ and $g_0$ are
in the bihomothetical relation (2.13) with
$q(s)=\displaystyle{\frac{p'(s)}{\alpha_0}}$.

Hence $Q(s)\subset N$ is an $\alpha$-Sasakian space form with
curvature tensor $K$ satisfying (2.16), where
$$\alpha=\displaystyle{\frac{p'}{p}>0,\quad
H=\frac{1}{p^2}\{H_0+3\alpha^2_0-3{p'}^2\}}.\leqno(3.26)$$

Replacing the tensor $K$ into (3.24) we find
$$R_{\vert \Delta}= \displaystyle{\frac{1}{4}\,
\left(H-\frac{{p'}^2}{p^2}\right)
(\pi_1+\pi_2)_{\vert\Delta}-\frac{1}{4}\,\left(H-\frac{{p'}^2}{p^2}+
4\frac{p''}{p}\right){\pi_3}_{\vert\Delta}},\leqno(3.27)$$
where $\pi_1, \pi_2$ and $\pi_3$ are the tensors (2.12)
with respect to the structure
$(\textsl{g},\varphi, \tilde \xi, \tilde\eta)$.

Taking into account (3.18), (3.27) and (3.23) we obtain
$$R=a\pi + b\Phi +c\Psi,$$
where
$$a=H-\frac{p'^2}{p^2}, \quad
b=-2\left(H-\frac{p'^2}{p^2}+4\frac{p''}{p}\right),
\quad c=H-\frac{p'^2}{p^2}+5\frac{p''}{p}-\frac{p'''}{p'}.\leqno{(3.28)}$$
According to \cite{GM2} $(N,\textsl{g},J)$ is a K\"ahler manifold
of quasi-constant holomorphic sectional curvatures.

For the inverse, let the warped product  K\"ahler manifold (3.15) be of
quasi-constant holomorphic sectional curvatures. Then any integral
submanifold $Q(s)$ of the distribution $\Delta$ is an
$\alpha$-Sasakian space form \cite{GM3}. Since $Q(s)$ and $Q_0$
are in bihomothetical correspondence (2.13), then the base manifold
$Q_0$ is an $\alpha_0$-Sasakian space form.
\hfill{\bf QED}
\vskip 2mm
Thus we constructed the following
\vskip 1mm
\noindent
{\bf Examples} of K\"ahler manifolds of quasi-constant holomorphic
sectional curvatures with $B_0$-distribution:
$$\begin{array}{l}
\vspace{2mm}
All \; warped \; product \; Kaehler\; manifolds \;
(N=Q_0\times {\R},\textsl{g}, J,D), \; whose \\
\vspace{2mm}
underlying \; base \; manifold \, \, Q_0 \, \, is \; an \;
\alpha_0{\rm-}Sasakian \; space \; form.
\end{array} \leqno{(E)}$$

\begin{rem}\label{R:3.2}
From (3.17), (3.25), (3.28) and (2.18) it follows that
$$a+k^2=H+3\alpha^2=\frac{H_0+3\alpha_0^2}{p^2}. $$
Hence, the type of the underlying $\alpha_0$-Sasakian space form
$Q_0$ (i.e. $H_0+3\alpha_0^2 \gtreqqless 0$) determines
uniquely the type of the corresponding
warped product K\"ahler manifold of quasi-constant holomorphic sectional
curvatures (i.e. $a+k^2 \gtreqqless 0$).
\end{rem}
Now we shall prove the main theorem in this section.
\begin{thm}\label{T:3.1}
Any K\"ahler manifold of quasi-constant holomorphic sectional
curvatures with $B_0$-distribution, satisfying one of the conditions
$$a+k^2>0, \quad a+k^2=0, \quad a+k^2<0,$$
locally has the structure of a warped product K\"ahler manifold from
the Examples (E).
\end{thm}

{\it Proof.} Let $(M,g,J,D)\,(\dim M=2n\geq 6)$ be a K\"ahler
$QCH$-manifold with $B_0$-distribution $D$ and functions $k,\,p^*,\,a,$
which satisfy one of the inequalities in the theorem and $m_0$ be an
arbitrary point in $M$. The Levi-Civita connection of the metric
$g$ is denoted as usual by $\nabla$.

Since $\eta$ is closed, we can find a coordinate neighborhood $U$ about
$m_0$ with coordinate functions $(w^1,...,w^{2n-1};s)$, where
$(w^1,...,w^{2n-1})$ are in a domain $W \in {\R}^{2n-1}, \;
s\in I'=(-\varepsilon, \varepsilon), \; \varepsilon >0$ and $\eta=ds$.
The integral submanifolds $Q(s)$ of $\Delta$ in $U$ are determined by
$Q(s): s=const \in I'$. Especially the submanifold $Q_0$ of $\Delta$
through $m_0$ in $U$ is given by $Q_0: s=0$. We have shown that
the integral submanifold $Q_0$ of $\Delta$ carries an
$\alpha_0$-Sasakian structure
$(g_0,\varphi,\tilde\xi_0,\tilde\eta_0)$ with
$\displaystyle{\alpha_0=\frac{k(0)}{2}}$ and constant
$\varphi$-holomorphic sectional curvatures
$\displaystyle{H_0=a(0)+\frac{k^2(0)}{4}}$. With the help of the local
base $\displaystyle{\left\{\frac{\partial}{\partial w^i}\right\}}, \;
i=1,...,2n-1$, of $Q_0$ we construct a frame field $\{e_{\beta},
\varphi \, e_{\beta},\tilde\xi_0\}, \; \beta=1,...,n-1$, which is a
$\varphi$-base for $T_mQ_0$ at every point $m\in Q_0$. We note
that this frame field is a $\varphi$-base on every $Q(s), \; s\in I'$.

We need the following lemma.
\begin{lem}\label{L:3.3}
Let $\{e_{\beta},\varphi \, e_{\beta},\tilde\xi_0\}, \; \beta=1,...,n-1$
be a frame field on $Q_0$, which is a $\varphi$-base at every point of
$Q_0$. There exist unique positive functions
$\lambda(s), \, \mu(s), \, \, s\in I'$, such that
$\lambda(0)=\mu(0)=1$ and the vector fields
$\tilde e_{\beta}=\lambda\, e_{\beta},\;\varphi \,
\tilde e_{\beta}=\lambda\, \varphi \,e_{\beta},\;
\tilde\xi=\mu \,\tilde \xi_0, \; \beta=1,...,n-1$, in ${\X}U$
are parallel along the integral curves of $\xi$.
\end{lem}
{\it Proof.} The conditions
$\left[\displaystyle{\frac{\partial}{\partial w^i},
\frac{\partial}{\partial s}}\right]=0$
imply that $[e_{\beta},\xi]=[\varphi \,e_{\beta},\xi]=
[\tilde \xi_0,\xi]=0, \; \beta= 1,...,n-1$, or equivalently
$$\nabla_{\xi}\,e_{\beta}=\frac{k}{2}\,e_{\beta},\quad
\nabla_{\xi}\,(\varphi \,e_{\beta})=\frac{k}{2}\,
(\varphi \,e_{\beta}),\quad
\nabla_{\xi}\,\tilde \xi_0=-p^* \tilde \xi_0.$$
The last equalities determine uniquely the functions
$$\lambda(s)=e^{\displaystyle{-\int_0^s\,\frac{k}{2}\,ds}},\quad
\mu(s)=e^{\displaystyle{\int_0^s\,p^*(s)\,ds}}, \quad s\in I',
\leqno{(3.29)}$$
such that
$$\nabla_{\xi}\,(\lambda \, e_{\beta})=\nabla_{\xi}\,\tilde e_{\beta}=0,
\quad \nabla_{\xi}\,\lambda (\varphi \,e_{\beta})=
\nabla_{\xi}\,(\varphi \,\tilde e_{\beta})=0,\quad
\nabla_{\xi}\,(\mu \,\tilde \xi_0)=\nabla_{\xi}\,\tilde \xi=0,$$
which proves the lemma.
\hfill{\bf QED}

We continue the proof of the theorem.

Because of (2.6) and (3.29) the functions $\lambda(s),\, \mu(s)$
are related as follows:
$$\frac{1}{\mu}=\frac{1}{\alpha_0}\,\frac{1}{\lambda}
\left(\frac{1}{\lambda}\right)'.\leqno{(3.30)}$$

Since ${\rm span}\left\{\displaystyle{\frac{\partial }
{\partial w^i}}\right\}=T_mQ(s)$ at any point
$m\in Q(s), \; s$ - fixed in $I'$, we can consider the metric
$g_0$ acting on $T_mQ(s)$ which is identified with the metric
$g_0$ of the basic leaf $Q_0$.

We give the relation between the metric $g(s)$ and the metric
$g_0$ in $U$.

Let $\tilde e_{\beta},\;\varphi \, \tilde e_{\beta},\;
\tilde\xi, \; \beta=1,...,n-1$, be as in Lemma \ref{L:3.3}. Since
$\tilde e_{\beta},\;\varphi \, \tilde e_{\beta},\;\beta=1,...,n-1$,
are parallel along the integral curves of $\xi$, we have
$$g(s)(\tilde e_{\beta},\tilde e_{\beta})=\lambda^2(s)
g(e_{\beta},e_{\beta})=1 =g_0(e_{\beta},e_{\beta}).$$
Hence
$$g(s)(e_{\beta},e_{\beta})=\frac{1}{\lambda^2}\,g_0(e_{\beta},e_{\beta}).$$
Similarly we have
$$g(s)(\varphi\,e_{\beta},\varphi\,e_{\beta})=\frac{1}{\lambda^2}\,
g_0(\varphi\,e_{\beta},\varphi\,e_{\beta}),$$
$$g(s)(\tilde\xi_0,\tilde\xi_0)=\frac{1}{\mu^2}\,
g_0(\tilde\xi_0,\tilde\xi_0).$$
Thus we obtain
$$g(s)=\frac{1}{\lambda^2}\,(g_0-\tilde\eta_0\otimes\tilde\eta_0)+
\frac{1}{\mu^2}\,\tilde\eta_0\otimes\tilde\eta_0.\leqno{(3.31)}$$
Putting $p(s)=\displaystyle{\frac{1}{\lambda(s)}}, \; s\in I'$, in
view of (3.30) and (3.31), we get
$$g(s)=p^2(g_0-\tilde\eta_0\otimes\tilde\eta_0)+
\left(\frac{pp'}{\alpha_0}\right)^2
\tilde\eta_0\otimes\tilde\eta_0.\leqno{(3.32)}$$

Further the proof of the theorem ends in the following scheme.

Let $Q_0\times I'$ be the standard product manifold, defined in the
coordinate set $W\times I'$. In this case we denote
$\displaystyle{\xi'=\frac{\partial}{\partial s}},\;
\eta'=ds,\;\Delta'={\rm span}\displaystyle{\left\{\frac{\partial}
{\partial w^1},...,\frac{\partial}{\partial w^{2n-1}}\right\}}$.
Then every leaf of $Q_0\times I'$ carries the $\alpha_0$-Sasakian
structure $(g_0,\varphi,\tilde\xi_0,\tilde\eta_0)$ of the base
manifold $Q_0$. Using the function $\displaystyle{p(s)=\frac{1}
{\lambda(s)}}$, where $\lambda(s)$ is given by (3.29),
we endow the manifold $U'=Q_0\times I'$ with the warped
product K\"ahler structure $(\textsl{g}',J',\xi',\eta')$ as in (3.15).

We shall show that the natural coordinate diffeomorphism
$$F: \; U \; \rightarrow \; U'$$
is an equivalence, i.e. $F$ preserves the structures $g,J,\xi$
(consequently $\eta, \Delta, D$).

By definition $F$ preserves the vector field
$\xi=\displaystyle{\frac{\partial}{\partial s}}$,
the 1-form $\eta=ds$ and the distribution $\Delta$.

From the condition that $F$ preserves the vector fields
$\displaystyle{\left\{\frac{\partial}{\partial w^i}\right\}},
\; i=1,...,2n-1$, it follows that $F$ preserves the $\varphi$-base
$\{e_{\beta},\varphi \, e_{\beta},\tilde\xi_0\}, \; \beta=1,...,n-1$
in ${\X}U$. This means that $F$ preserves the structure $\varphi$
on every leaf $Q(s), \; s\in I'$. Taking into account (3.9),
we conclude that $F$ preserves the complex structure $J$.

Finally, the metric $g(s)$ on $Q(s)\subset U, \; s\in I'$ has the form
(3.32). On the other hand, according to (3.15), the restriction of
the metric $\textsl{g}'$ of $U'$ onto the leaves $Q'(s)$ is given
by the same formula (3.32), written with respect to
the metric $g_0$ of the base $Q_0$. Hence $F$ is also an isometry.
\hfill{\bf QED}

\section{On the local theory of Bochner-K\"ahler manifolds}
\vskip 2mm
Let $(M,g,J) \; (\dim \,M = 2n \geq 4)$ be a K\"ahler manifold.
In the next calculations we shall use the complexification
$T_m^{\C}M,  \, m\in M$ and its standard splitting
$$T_m^{\C}M = T_m^{1,0}M \oplus T_m^{0,1}M.$$
Any complex basis of $T_m^{1,0}M$ will be denoted by
$\{Z_{\alpha}\}, \, \alpha = 1,...,n$, and the conjugate basis
$\{Z_{\bar\alpha} = \overline{Z_{\alpha}}\}, \, \bar\alpha =
\bar 1,...,\bar n$, will span $T_m^{0,1}M$. Unless otherwise stated,
the Greek indices $\alpha, \beta, \gamma, \delta, \varepsilon$ will
run through $1,...,n$.

We recall that the Bochner curvature operator $B$ acts on the
curvature tensor $R$ of $(M,g,J)$ with respect to a complex base
as follows:
$$\begin{array}{ll}
(B(R))_{\alpha\bar\beta\gamma\bar\delta}=
&R_{\alpha\bar\beta\gamma\bar\delta}-
\displaystyle{
\frac{1}{n+2}\,(g_{\alpha\bar\beta}\rho_{\gamma\bar\delta}
+g_{\gamma\bar\beta}\rho_{\alpha\bar\delta} +
g_{\gamma\bar\delta}\rho_{\alpha\bar\beta}
+ g_{\alpha\bar\delta}\rho_{\gamma\bar\beta})}\\
[3mm]& \displaystyle{+
\frac{\tau}{2(n+1)(n+2)}\,(g_{\alpha\bar\beta}g_{\gamma\bar\delta}+
g_{\gamma\bar\beta}g_{\alpha\bar\delta})}.
\end{array}$$

The manifold $(M,g,J)$ is said to be Bochner-K\"ahler if its
Bochner curvature tensor $B(R)$ vanishes identically, which is
equivalent to the curvature identity:
$$\begin{array}{ll}
R_{\alpha\bar\beta\gamma\bar\delta} =&\displaystyle{ \frac{1}{n+2}
\,(g_{\alpha\bar\beta}\rho_{\gamma\bar\delta}
+g_{\gamma\bar\beta}\rho_{\alpha\bar\delta} +
g_{\gamma\bar\delta}\rho_{\alpha\bar\beta}
+ g_{\alpha\bar\delta}\rho_{\gamma\bar\beta})}\\
[3mm]& \displaystyle{-
\frac{\tau}{2(n+1)(n+2)}\,(g_{\alpha\bar\beta}g_{\gamma\bar\delta}+
g_{\gamma\bar\beta}g_{\alpha\bar\delta})}.
\end{array}\leqno{(4.1)}$$

The curvature tensor of any K\"ahler manifold as a consequence of
the second Bianchi identity satisfies the equalities:
$$\nabla_{\alpha}\,\rho_{\gamma\bar\beta} =
\nabla_{\gamma}\,\rho_{\alpha\bar\beta}=
\nabla_{\varepsilon}\,R_{\alpha\bar\beta\gamma}^{\varepsilon}.$$
The other meaning of the identity \;
$\nabla_{\alpha}\,\rho_{\gamma\bar\beta} =
\nabla_{\gamma}\,\rho_{\alpha\bar\beta}$ \;
is that the Ricci form $\rho(JX,Y), \; X,Y \in {\X}M$ is closed, i. e.
$$(\nabla_X\,\rho)(JY,Z)+(\nabla_Y\,\rho)(JZ,X)+(\nabla_Z\,\rho)(JX,Y)=0,
\quad X,Y,Z \in {\X}M. \leqno{(4.2)}$$
On the other hand (4.1) implies
$$\nabla_{\varepsilon}\,R_{\alpha\bar\beta\gamma}^{\varepsilon}=
\frac{2}{n+2}\,\nabla_{\alpha}\,\rho_{\gamma\bar\beta}+
\frac{n}{2(n+1)(n+2)}\,
(\tau_{\alpha}g_{\gamma\bar\beta} +\tau_{\gamma}g_{\alpha\bar\beta}).$$
Hence
$$\nabla_{\alpha}\,\rho_{\gamma\bar\beta}= \frac{1}{2(n+1)}\,
(\tau_{\alpha}g_{\gamma\bar\beta} +
\tau_{\gamma}g_{\alpha\bar\beta}).\leqno{(4.3)}$$
\begin{rem}
The identity (4.3) in view of (4.2) is equivalent to
$$\begin{array}{ll}
(\nabla_X\,\rho)(Y,Z)=&\displaystyle{\frac{1}{4(n+1)}\,\{2d\tau(X)g(Y,Z)
+d\tau(Y)g(X,Z)+d\tau(Z)g(X,Y)}\\
[3mm]&+d\tau(JY)g(X,JZ)+d\tau(JZ)g(X,JY)\}, \quad X,Y,Z \in {\X}M.
\end{array}$$
\end{rem}

From (4.1) and (4.3) it follows that the covariant derivative of the
curvature tensor $R$ of any Bochner-K\"ahler manifold satisfies the
following identity:
$$\nabla_{\alpha}\,R_{\beta \bar\varepsilon \gamma \bar\delta}=
\frac{1}{(n+1)(n+2)}\,(\tau_{\alpha}\pi_{\beta \bar \varepsilon
\gamma \bar\delta} +
\tau_{\beta}\pi_{\alpha\bar\varepsilon\gamma\bar\delta}+
\tau_{\gamma}\pi_{\beta
\bar\varepsilon\alpha\bar\delta}).\leqno{(4.4)}$$

Studying the integrability conditions for (4.3) we obtain the
basic formulas connecting derivatives of the scalar curvature
$\tau$ with curvature properties of the given manifold.
\begin{prop}\label{P:4.1}
Let $(M,g,J)$ be a Bochner-K\"ahler manifold and $T=grad \, \tau$.
Then

(i) the vector field\, $\displaystyle {\frac{T-iJT}{2}}$\, is holomorphic;

(ii) $\displaystyle{(n+2)\nabla_{\alpha}\,\tau_{\bar\beta}
+2(n+1)\rho_{\alpha\bar\beta}^2 -
\tau \rho_{\alpha\bar\beta}=\frac{(n+2)\Delta \tau
+ 2(n+1)\Vert \rho \Vert^2 - \tau^2}
{2n}\, g_{\alpha\bar\beta}}$,

\hspace{5mm}
where  $\rho_{\alpha \bar\beta}^2=
\rho_{\alpha}^{\varepsilon}\rho_{\varepsilon\bar\beta}$;
\vspace {1mm}

(iii) $2\rho(X,T)=-X(\Delta\tau), \quad X \in {\X}M.$
\end{prop}

{\it Proof.} Applying the standard Ricci formula in Riemannian geometry
$$\nabla_i\nabla_j\,\rho_{km}-\nabla_j\nabla_i\,\rho_{km}
=-R_{ijk}^s \rho_{sm}-R_{ijm}^s \rho_{ks}$$
and taking into account (4.3) we obtain
$$\nabla_{\alpha}\nabla_{\beta}\,\rho_{\gamma\bar\delta}
- \nabla_{\beta}\nabla_{\alpha}\,\rho_{\gamma\bar\delta} =
\displaystyle{
\frac{1}{2(n+1)}\,(g_{\beta\bar\delta}\nabla_{\alpha}\,\tau_{\gamma}-
g_{\alpha\bar\delta}\nabla_{\beta}\,\tau_{\gamma})}=0,
\leqno{(4.5)}$$
$$\begin{array}{lll}
\nabla_{\alpha}\nabla_{\bar\beta}\,\rho_{\gamma\bar\delta}-
\nabla_{\bar\beta}\nabla_{\alpha}\,\rho_{\gamma\bar\delta} &=&
\displaystyle{-\frac{1}{n+2}\,
(g_{\gamma\bar\beta}\rho^2_{\alpha\bar\delta}-
g_{\alpha\bar\delta}\rho^2_{\gamma\bar\beta})}\\
[2mm] & &+
\displaystyle{
\frac{\tau}{2(n+1)(n+2)}\,(g_{\gamma\bar\beta}\rho_{\alpha\bar\delta}-
g_{\alpha\bar\delta}\rho_{\gamma\bar\beta})}\\
[3mm]&
=&\displaystyle{
\frac{1}{2(n+1)}\,(g_{\gamma\bar\beta}\nabla_{\alpha}\,\tau_{\bar\delta}-
g_{\alpha\bar\delta}\nabla_{\gamma}\,\tau_{\bar\beta})}.
\end{array} \leqno{(4.6)}$$
After taking a trace in (4.5) we find
$$\nabla_{\beta}\,\tau_{\gamma} = 0. \leqno{(4.7)}$$
This equality implies that the $(1,0)$-part\,
 $\displaystyle{\frac{T-iJT}{2}}$\,
of the vector field $T = grad \, \tau$ is holomorphic, which proves (i).

In a similar way (4.6) implies
$$(n+2)\nabla_{\gamma}\,\tau_{\bar\beta}+2(n+1)\rho_{\gamma\bar\beta}^2 -
\tau \rho_{\gamma\bar\beta}=\frac{(n+2)\Delta \tau + 2(n+1)\Vert
\rho \Vert^2 - \tau^2} {2n}\, g_{\gamma\bar\beta}, \leqno{(4.8)}$$
which is (ii).

To prove (iii) we consider the identity
$$\nabla_{\bar\alpha}\nabla_{\beta}\,\tau_{\gamma}-
\nabla_{\beta}\nabla_{\bar\alpha}\,\tau_{\gamma}=
-R_{\bar\alpha\beta\gamma}^{\varepsilon}\tau_{\varepsilon}.$$
Keeping in mind (4.7) we take a trace in the last equality and get
$$-\frac{1}{2}(\Delta \tau)_{\beta} =
\rho_{\beta}^{\varepsilon}\tau_{\varepsilon},$$
which is
$$2\rho(X,T)=-X(\Delta \tau), \quad X \in {\X}M.$$
\hfill{\bf QED}
\begin{rem}
The equalities
$$\nabla_{\alpha}\,\tau_{\beta}=0, \quad \nabla_{\alpha}\,\tau_{\bar\beta}-
\nabla_{\bar\beta}\,\tau_{\alpha}=0$$
imply $JT$ is an analytic and Killing vector field. Hence $JT$ generates
a local one-parameter group of local holomorphic motions.
\end{rem}

\begin{prop}\label{P:4.2}
On any Bochner-K\"ahler manifold the function
$$\Vert \rho \Vert^2 -\frac{\tau^2}{2(n+1)} + \frac{\Delta \tau}{n+1}$$
is a constant.
\end{prop}
{\it Proof.} From the equality (4.3) we obtain successively
$$\Vert \rho \Vert^2_{\alpha}
= 4\rho^{\beta\bar\gamma}\nabla_{\alpha}\,\rho_{\beta\bar\gamma}
= \frac{1}{n+1}\,(\frac{\tau^2}{2}-\Delta \tau)_{\alpha}.$$
Hence
$$\Vert \rho \Vert^2 - \frac{\tau^2}{2(n+1)} + \frac{\Delta \tau}{n+1}
= const.$$
\hfill{\bf QED}

We set
$$\Vert \rho \Vert^2 - \frac{\tau^2}{2(n+1)} + \frac{\Delta \tau}{n+1}
= {\B}\leqno{(4.9)}$$
and call ${\B}$ {\it the Bochner constant} of the manifold.
\vskip 2mm
It is clear that the 1-form $d\tau$ is of basic importance in
the above formulas. The equality (4.3) shows that the conditions
$\tau = const$ and $\nabla\rho = 0$ are equivalent on a
Bochner-K\"ahler manifold. Because of the structural theorem in
\cite {TL} the case $B(R)=0, \; d\tau =0$ can be considered as
well-studied.

According to \cite{K,B} it follows that any compact Bochner-K\"ahler
manifold satisfies the condition $d\tau =0 \; (\nabla R=0)$.

Further we consider the general case of Bochner-K\"ahler manifolds
$$d\tau \neq 0 \quad \text{for all points} \; p \in M.\leqno{(4.10)}$$
This condition allows us to introduce the frame field
$$\left\{ \xi = \frac{grad \, \tau}{\Vert d\tau \Vert}
= \frac{T}{\Vert d\tau \Vert}, \quad
J\xi = \frac{Jgrad \, \tau}{\Vert d\tau \Vert} = \frac{JT}{\Vert
d\tau \Vert}\right\}$$
and consider the $J$-invariant distributions $D_{\tau}$ and
$D^{\perp}_{\tau}=span\{\xi,J\xi\}$.

Thus our approach to the local theory of Bochner-K\"ahler manifolds is
to investigate them as K\"ahler manifolds $(M,g,J,D_{\tau})$ endowed with
a $J$-invariant distribution $D_{\tau}$ generated by the K\"ahler
structure $(g,J)$. In what follows we call this distribution
{\it the scalar distribution} of the manifold.

The scalar distribution $D_{\tau}$ of any Bochner-K\"ahler manifold
carries the functions
$$\varkappa:=R(\xi,J\xi,J\xi,\xi), \quad
\sigma:=\rho(\xi,\xi)=\rho(J\xi,J\xi).$$
Then $\varkappa, \,\sigma$ and $\tau$ determine the functions
$$ \begin{array}{l}
\displaystyle{ a = \frac{\tau - 4 \sigma + 2\varkappa}{n(n-1)},
\quad b = \frac{- 2\tau +4(n+2)\sigma - 4(n+1)\varkappa}{n(n-1)}},\\
[4mm]
\displaystyle{c = \frac{\tau - 4(n+1)\sigma +
(n+1)(n+2)\varkappa}{n(n-1)}}.
\end{array}$$

Calculating \,$\varkappa$\, from (4.1) we obtain
$$0=\tau-4(n+1)\sigma+(n+1)(n+2)\varkappa=n(n-1)\,c.
\leqno{(4.11)}$$

Below we establish some properties of the distributions $D_{\tau},\,
D^{\perp}_{\tau}$ and $\Delta_{\tau} \; (\perp \xi)$.

Let $\{Z_{\lambda}\}, \, \lambda = 1,...,n-1$ be a basis for
$D^{1,0}(m)$. The basis $\{Z_0,Z_{\lambda}\}, \,\lambda = 1,...,n-1$,
where $\displaystyle{Z_0 = \frac{\xi -iJ\xi}{2}}$, is said to be a
{\it special complex basis} for $T_m^{1,0}M$. Then
$\{Z_{\bar 0}, Z_{\bar\lambda}\}, \,
\bar\lambda = \bar 1,...,\overline{n-1}$, is a special complex
basis for $T_m^{0,1}M$. The Greek indices $\lambda, \mu, \nu,
\varkappa, \sigma$ will run through $1,...,n-1.$

With respect to special complex bases the 1-forms $\eta$ and
$\tilde\eta$ have the following components:
$$\eta_{\alpha} = g_{\alpha \bar 0}, \quad \eta_{\bar\alpha} =
g_{\bar\alpha 0}, \quad \tilde\eta_{\alpha} = -i\eta_{\alpha},
\quad \tilde\eta_{\bar\alpha} = i\eta_{\bar\alpha},$$
$$\eta_{\lambda} = \eta_{\bar\lambda} = 0, \quad \eta_0 =
\eta_{\bar 0} = g_{0 \bar 0} = \frac{1}{2}. $$

We introduce the following functions and 1-forms associated with
the vector fields $\nabla_{\xi}\,\xi$ and $\nabla_{J\xi}\,J\xi$:
$$p = g(\nabla_{\xi}\,\xi,J\xi), \quad p^* = g(\nabla_{J\xi}\,J\xi,\xi),
\leqno{(4.12)}$$
$$\theta(X) = g(\nabla_{\xi}\,\xi,X) - p\,\tilde\eta(X), \quad
\theta^*(X) = g(\nabla_{J\xi}\,J\xi,X) - p^*\,\eta(X), \quad
X \in T_m M. \leqno{(4.13)}$$
It is clear that $\theta(X) = \theta(x_0), \; \theta^*(X) =
\theta^*(x_0)$, where $x_0 = X - \tilde\eta(X)J\xi - \eta(X)\xi.$

Taking into account that $Z_0 = \displaystyle{\frac{\xi
-iJ\xi}{2}, \; Z_{\bar 0} = \frac{\xi + iJ\xi}{2}}$ we find
$$\nabla_0\,\eta_{\lambda} = \frac{1}{2}\,(\theta_{\lambda} +
\theta^*_{\lambda}), \quad \nabla_{\bar 0}\,\eta_{\lambda} =
\frac{1}{2}\,(\theta_{\lambda} - \theta^*_{\lambda}), \leqno{(4.14)}$$
$$\nabla_0\,\eta_0 = \frac{p^* - ip}{4}, \quad \nabla_0\,\eta_{\bar 0} =
\frac{-p^* + ip}{4}.\leqno{(4.15)}$$

Since the distribution $\Delta_{\tau}$ is involutive then
$$\nabla_{\lambda}\,\eta_0=-\nabla_{\lambda}\,\eta_{\bar0}
= \frac{1}{2}\,\theta^*_{\lambda}.\leqno{(4.16)}$$
From $d\tau= \Vert d\tau \Vert \, \eta$ we have
$$(\nabla_X\,d\tau)Y=\Vert d\tau \Vert (\nabla_X\,\eta)Y+
X(\Vert d\tau \Vert)\eta(Y), \quad X,Y\in{\X}M.\leqno{(4.17)}$$
Then (4.7) and (4.17) imply that
$$\Vert d\tau \Vert_{\alpha}\eta_{\beta}+
\Vert d\tau \Vert \nabla_{\alpha}\,\eta_{\beta} =0.$$
The last equality, (4.10), (4.14), (4.16) and (4.15) give
$$\begin{array}{l}
\nabla_{\lambda}\,\eta_{\mu}= 0;\\
[2mm]
\theta_{\lambda}+\theta^*_{\lambda} = 0, \quad
\theta_{\lambda} = (\ln \Vert d\tau \Vert)_{\lambda};\\
[2mm]
p = -J\xi(\ln \Vert d\tau \Vert), \quad p^* = -\xi (\ln \Vert
d\tau \Vert).
\end{array} \leqno{(4.18)}$$

In a similar way the equality $\nabla_{\alpha}\,\tau_{\bar\beta}
=\nabla_{\bar\beta}\,\tau_{\alpha}$ implies
$$\begin{array}{l}
\nabla_{\lambda}\,\eta_{\bar\mu} = \nabla_{\bar\mu}\,\eta_{\lambda};\\
[2mm]
p = J\xi(\ln \Vert d\tau \Vert). \end{array}\leqno{(4.19)}$$

From (4.17) we calculate $div\,\xi=
\displaystyle{\frac{\Delta \tau}{\Vert d\tau \Vert}}+p^*.$
On the other hand $div_0\,\xi=div\,\xi-p^*$. Hence
$$div_0\,\xi=\frac{\Delta \tau}{\Vert d\tau \Vert}.$$

Summarizing the conditions (4.18) and (4.19) we get

\begin{prop}\label{P:4.3}
The scalar distribution of any Bochner-K\"ahler manifold has the
following properties:
$$\begin{array}{l}
i) \; \nabla_{\lambda}\,\eta_{\mu} = 0, \quad
\nabla_{\lambda}\,\eta_{\bar\mu} = \nabla_{\bar\mu}\,\eta_{\lambda},\\
[2mm]
ii) \; \theta_{\lambda} +\theta^*_{\lambda} = 0, \quad \theta_{\lambda}=
(\ln \Vert d\tau \Vert)_{\lambda},\\
[2mm]
iii) \; p=0=J\xi(\ln \Vert d\tau \Vert), \quad p^* =
-\xi(\ln \Vert d\tau \Vert),\\
[2mm]
iv) \; div_0\,\xi=\displaystyle{\frac{\Delta \tau}{\Vert d\tau \Vert}}.
\end{array}\leqno{(4.20)}$$
\end{prop}

Thus we have
\vskip 2mm
\noindent
{\it The local theory of any Bochner-K\"ahler manifold with scalar
distribution $D_{\tau}$ is determined by:
the symmetric tensor $\nabla_{\lambda}\eta_{\bar\mu}$,
the 1-form $\theta_{\lambda}=(\ln \Vert d\tau \Vert)_{\lambda}$
and the function $p^* =-\xi(\ln \Vert d\tau \Vert).$}
\vskip 2mm
Any additional conditions for the above mentioned objects give
rise to special classes of Bochner-K\"ahler manifolds.

\section{Bochner-K\"ahler manifolds whose scalar distribution is
a $B_0$-distribution}

The first step in our study of Bochner-K\"ahler manifolds with
respect to their scalar distribution is to study the basic class
of Bochner-K\"ahler manifolds whose scalar distribution is a
$B_0$-distribution.We start with the following curvature
characterization of these manifolds.

\begin {prop}\label{P:5.1}
Let $(M,g,J) \; (\dim M = 2n \geq 6)$ be a Bochner-K\"ahler manifold
with $d\tau \neq 0$. If the scalar distribution $D_{\tau}$ is a
$B_0$-distribution, then the manifold is of quasi-constant holomorphic
sectional curvatures and
$$R = a\pi + b\Phi, \quad b \neq 0.$$
\end{prop}

{\it Proof.} Let $D_{\tau}$ be a $B_0$-distribution,
i.e. we have to add the conditions
$$\begin{array}{l}
\displaystyle{\nabla_{\lambda}\,\eta_{\bar\mu} =
\frac{k}{2}\,g_{\lambda\bar\mu}, \quad
k = \frac{div_0 \, \xi }{n-1} \neq 0,}\\
[3mm]
\theta_{\lambda}= (\ln \Vert d\tau \Vert)_{\lambda} = 0
\end{array}\leqno{(5.1)}$$
to the equations (4.20).

Then (2.5) implies that:
$$\begin{array}{ll}
R(X,Y)\xi = & \displaystyle{-\frac{1}{4}\,(k^2+2kp^*)\{\eta(X)Y-\eta(Y)X}\\
[3mm]
& \displaystyle{-\tilde\eta(X)JY+\tilde\eta(Y)JX+2g(JX,Y)J\xi\}}\\
[2mm]
& \displaystyle{-\frac{1}{2k}\,\xi(k^2+2kp^*)
\{\eta(X)\tilde\eta(Y)-\eta(Y)\tilde\eta(X)\}J\xi},
\end{array}\leqno{(5.2)}$$
$X,Y \in {\X}M$.

Taking a trace in (5.2) we find
$$\rho(X,\xi)=\left\{\frac{1}{2k}\,\xi(k^2+2kp^*)+
\frac{n+1}{2}\,(k^2+2kp^*)\right\}\eta(X),
\quad X \in {\X}M.\leqno{(5.3)}$$
From (5.2) and (5.3) it follows that
$$\varkappa = R(\xi,J\xi,J\xi,\xi)=
\frac{1}{2k}\,\xi(k^2+2kp^*)+(k^2+2kp^*),\leqno{(5.4)}$$
$$\sigma = \rho(\xi,\xi)
=\frac{1}{2k}\,\xi(k^2+2kp^*)+\frac{n+1}{2}\,(k^2+2kp^*).\leqno{(5.5)}$$
The equalities (4.1) and (5.2) give two expressions for the component\,
$R_{0\bar\mu\nu\bar 0}$ \,of the curvature tensor $R$. Comparing these
expressions and taking into account (4.11) we obtain
$$\rho_{\nu\bar\mu}= \frac{\tau - 2\sigma}{2(n-1)}\,g_{\nu\bar\mu}.
\leqno{(5.6)}$$
The conditions (5.3) and (5.6) imply
$$\rho = \frac{\tau - 2\sigma}{2(n-1)}\,g+\frac{2n\sigma - \tau}{2(n-1)}\,
(\eta\otimes\eta+\tilde\eta\otimes\tilde\eta),\leqno{(5.7)}$$
which combined with (4.1) gives
$$R = a\pi + b\Phi.\leqno{(5.8)}$$
Thus we obtained that $(M,g,J,D_{\tau})$ is a K\"ahler $QCH$-manifold with
$B_0$-distribution. Applying the second Bianchi identity to (5.8)
it follows that (cf Theorem 3.5 in \cite{GM2})
$$da=\frac{1}{2}\,b\,k\,\eta, \quad db=b\,k\,\eta. \leqno{(5.9)}$$
From (5.8) we calculate
$$\tau=(n+1)(na+b).\leqno{(5.10)}$$
On the other hand, it follows from (5.9) that
$$d(2a-b)=0.\leqno{(5.11)}$$
Combining (5.10) and (5.11) we get
$$d\tau = \frac{1}{2}\,(n+1)(n+2)db\neq 0,\leqno{(5.12)}$$
which implies that $b\neq0$.
\hfill{\bf QED}
\vskip 1mm
Taking into account (5.11), we denote
$${\bb}_0=\frac{2a-b}{2}=const.\leqno{(5.13)}$$
\vskip 2mm
Thus {\it any Bochner-K\"ahler manifold whose scalar distribution is
a $B_0$-distribution has two geometric constants: the Bochner constant
${\B}$ given by $(4.9)$ and ${\bb}_0$ given by $(5.13)$.}
\vskip 2mm
Next we find expressions for the functions $\varkappa,\; \sigma,\;
a,\; b$ and $k$ by $\tau$.
\begin{lem}\label{L:5.1}
Let $(M,g,J)\;(\dim M = 2n \geq 6)$ be a Bochner-K\"ahler manifold whose
scalar distribution is a $B_0$-distribution. Then
$$\varkappa=\frac{3\tau}{(n+1)(n+2)}-\frac{2(n-1){\bb}_0}{n+2}, \quad
\sigma=\frac{\tau}{n+1}-\frac{(n-1){\bb}_0}{2}, \leqno{(5.14)}$$
$$a=\frac{\tau}{(n+1)(n+2)}+\frac{2{\bb}_0}{n+2},\quad
b=\frac{2\tau}{(n+1)(n+2)}-\frac{2n{\bb}_0}{n+2},\leqno{(5.15)}$$
$$k^2=\frac{(n+1)^2(n+2)^2{\K}-
\{2\tau-(n+1)(n-2){\bb}_0\}^2}{4(n+1)(n+2)\{\tau-n(n+1){\bb}_0\}},
\leqno{(5.16)}$$
where
$${\K}= \frac{4{\B}+ n^2(2n+1){\bb}_0^2}{n(n+2)}. \leqno{(5.17)}$$
\end{lem}
{\it Proof.} Taking into account (2.3), (4.11) and (5.13) we find
(5.14) and (5.15).

From (5.2), (5.4) and (5.5) we get
$$R(x_0,\xi,\xi,x_0)= \frac{\sigma-\varkappa}{2(n-1)}
=\frac{1}{4}\,(k^2+2kp^*), \quad x_0 \in D,\quad g(x_0,x_0) = 1.$$
On the other hand, from (5.8) we have
$$R(x_0,\xi,\xi,x_0)=\frac{2a+b}{8}=
\frac{\sigma - \varkappa}{2(n-1)}, \quad x_0 \in D, \quad g(x_0,x_0)= 1.$$
Comparing the above equalities in view of (5.13) and (2.6) we find
$$\xi(k)= -\frac{1}{2}\,(k^2+b + {\bb}_0).\leqno{(5.18)}$$
Then (5.9) implies that
$$\frac{dk^2}{db}=\frac{2k\,\xi(k)}{\xi(b)}= -\frac{k^2+b+{\bb}_0}{b}.$$
The general solution of the last equation is
$$k^2=\frac{1}{2b}\,\{{\K}-(b+{\bb}_0)^2\}\leqno{(5.19)}$$
where ${\K}= const.$

We shall express the constant \,${\K}$\, by means of the constants\,
${\B}$ \,and \,${\bb}_0$\, of the manifold under consideration.

By using (5.12) and (5.9) we have
$$\Vert d\tau \Vert = \xi(\tau)=\frac{(n+1)(n+2)}{2}\,\xi(b)=
\frac{(n+1)(n+2)}{2}\,kb.\leqno{(5.20)}$$

Differentiating the equality
$\tau_{\alpha}=\Vert d\tau \Vert \eta_{\alpha}$ we find
$$\nabla_{\alpha}\,\tau_{\bar\beta}=
\Vert d\tau \Vert \nabla_{\alpha}\,\eta_{\bar\beta}
+ \Vert d\tau \Vert_{\alpha}\eta_{\bar\beta} \leqno{(5.21)}$$
with respect to a complex basis. Taking a trace in (5.21)
because of (5.20), (5.9) and (5.1) we get
$$\Delta \tau = \frac {(n+1)(n+2)}{2}\,b\,\{2\xi(k)+ (n+1)k^2\}.
\leqno{(5.22)}$$
Further we replace\, $\xi(k)$ \,from (5.18) and \,$k^2$\, from
(5.19) into (5.22) and obtain
$$\frac{\Delta \tau}{n+1}= \frac{n+2}{4}\,\{n{\K}-(n+2)b^2-
2(n+1)b{\bb}_0-n{\bb}_0^2\}.
\leqno{(5.23)}$$
From (5.10) and (5.13) it follows that
$$\frac{\tau^2}{2(n+1)}= \frac{n+1}{8}\,\{(n+2)^2b^2+4n(n+2)b{\bb}_0+
4n^2{\bb}_0^2\}.
\leqno{(5.24)}$$

In order to calculate the constant ${\B}$ it remains to find
$\Vert \rho \Vert ^2$. By using the equality
$$\rho_{\alpha\bar\beta}=\frac{(n+2)b +
2(n+1){\bb}_0}{4}\,g_{\alpha\bar\beta}+
\frac{n+2}{2}\,b\,\eta_{\alpha}\eta_{\bar\beta}$$
we get
$$\Vert \rho \Vert^2=\frac{(n+2)^2(n+3)}{8}\,b^2 +
\frac{(n+1)^2(n+2)}{2}\,b{\bb}_0+ \frac{n(n+1)}{2}\,{\bb}_0^2.
\leqno{(5.25)}$$

Counting (5.23), (5.24) and (5.25) we find
$${\B}= \Vert \rho \Vert^2-\frac{\tau^2}{2(n+1)}+ \frac{\Delta \tau}{n+1}=
\frac{n(n+2)}{4}\,{\K}-\frac{n^2(2n+1)}{4}\,{\bb}_0^2.$$
Hence
$${\K}= \frac{4{\B}+ n^2(2n+1){\bb}_0^2}{n(n+2)}. $$

Now (5.19), (5.15) and (5.17) imply the assertion. \hfill{\bf QED}
\vskip 1mm
It follows from (5.17) that the pair $({\B},{\bb}_0)$
determines the pair $({\K},{\bb}_0)$ and vice versa. Further in our
considerations we shall use the pair of geometric constants
\,$({\K},{\bb}_0)$.

From (5.13) and (5.19) we get
$$a+k^2=\frac{(n+1)(n+2)({\K}-{\bb}_0^2)}{4\{\tau-n(n+1){\bb}_0\}}.
\leqno{(5.26)}$$

In \cite{GM3} we have shown that a K\"ahler manifold with curvature
tensor satisfying (2.2) is Bochner flat if and only if $c=0$.

Another characterization is given by the following

\begin{lem}\label{L:5.2}
Let $(M,g,J,D) \; (\dim M = 2n \geq 6)$ be a K\"ahler $QCH$-manifold
with $B_0$-distribution $D$. Then the following conditions are equivalent:

i)\;\, $B(R)=0$,

ii)\; $2a-b=const.$
\end{lem}
{\it Proof.} Under the conditions of the lemma the curvature tensor \,$R$
\,has the form (2.2). Then the second Bianchi identity applied to \,$R$\,
(Theorem 3.5 \cite{GM2}) implies that
$$d(2a-b)=-4kc\,\eta.$$
Now the statement of the lemma follows immediately.
\hfill{\bf QED}
\vskip 2mm

\section{Warped product Bochner-K\"ahler metrics}

From Theorem \ref{T:3.1} and Proposition \ref{P:5.1} it follows that
the local theory of Bochner-K\"ahler manifolds whose scalar distribution
is a $B_0$-distribution is equivalent to the local theory of warped
product Bochner-K\"ahler manifolds.

In this section we give a complete description of the warped product
K\"ahler metrics of quasi-constant holomorphic sectional curvatures,
which are Bochner flat.

Let $(N,g,J,\xi) \; (\dim M = 2n \geq 6)$ be a warped product K\"ahler
manifold given as in (3.15).

Next we shall consider the metric (3.15) with respect to the function
$p(t), \; t\in I \subset {\R}$, satisfying the initial conditions
$p(t_0)=1, \; \displaystyle{\frac{dp}{dt}(t_0)=\alpha_0}$, in the form
$$\textsl{g}=p^2(t)\displaystyle{\left\{g_0+
\left(\frac{1}{\alpha_0}\frac{dp}{dt}-1\right)\,
\tilde\eta_0\otimes\tilde\eta_0\right\}+
\eta\otimes\eta}.\leqno{(6.1)}$$

Thus any warped product K\"ahler metric of quasi-constant holomorphic sectional curvatures
 is determined uniquely by the function $p(t)$ and the underlying
 $\alpha_0$-Sasakian space form $Q_0$.

Now we shall find the functions $p(t)$ for which the metric (6.1)
is Bochner-K\"ahler.

The functions $a(s)$ and $b(s)$ are given by (3.28). According to
Lemma \ref {L:5.2} the manifold $(N,g,J,\xi)$ is Bochner flat if and
only if $2a-b=const=2{\bb}_0$. We denote
$$d_0:=H_0+3\alpha_0^2.$$

Taking into account (5.18), (5.19), (3.17) and (5.13) we find
$$\left(\frac{p'}{p}\right)'-\left(\frac{p'}{p}\right)^2=
\frac{{\bb}_0}{4}-\frac{d_0}{2p^2}\leqno{(6.2)}$$
and
$$\left\{\left(\frac{p'}{p}\right)'-\left(\frac{p'}{p}\right)^2\right\}^2=
\frac{{\K}}{16}+\frac{d_0}{p^2}\left(\frac{p'}{p}\right)'.\leqno{(6.3)}$$
\vskip 2mm
\centerline {\it Warped product Bochner-K\"ahler metrics with \;
$a+k^2\neq 0 \; \; (d_0=H_0+3\alpha_0^2\neq 0).$}
\vskip 2mm
The equalities (6.2) and (6.3) imply
$$16d_0\,{p'}^2=({\bb}_0^2-{\K})\,p^4-4{\bb}_0d_0\,p^2+4d_0^2.
 \leqno{(6.4)}$$
Next we treat the function $p(s),\; s\in I_0$ with respect to the
parameter $t \in I$. Keeping in mind Remark \ref{R:3.1} we have
$$p'=\sqrt{\alpha_0\,\frac{dp}{dt}}.$$
Then (6.4) becomes
$$16\alpha_0\,d_0\frac{dp}{dt}=({\bb}_0^2-{\K})\,p^4-
4{\bb}_0d_0\,p^2+4d_0^2.
\leqno{(6.5)}$$

The initial conditions $p(t_0)=1$ and $\displaystyle{\frac{dp}{dt}(t_0)=
\alpha_0}$ \,give the following relation between the geometric
constants of the manifold:
$$16\alpha_0^2\,d_0+{\K}=({\bb}_0-2d_0)^2.\leqno{(6.6)}$$
Then from (6.5) we have
$$t=16\alpha_0\,d_0\int\frac{dp}{({\bb}_0^2-{\K})\,p^4-
4{\bb}_0d_0\,p^2+4d_0^2}.
\leqno{(6.7)}$$

Analyzing (6.7) and taking into account (6.6) we obtain the
following solutions of the equation (6.5):
\vskip 2mm
\centerline{\it I. The case \;$a+k^2>0$.}
\vskip 2mm
This case is determined by $d_0>0$, which is equivalent to
$H_0>-3\alpha_0^2$.
\vskip 2mm
Type 1.\quad ${\K}>0, \quad -\sqrt {\K}<{\bb}_0<\sqrt {\K} \quad
({\bb}_0^2-{\K}<0)$.

$$t=\frac{2\sqrt 2\,\alpha_0\sqrt{\sqrt {\K}-{\bb}_0}}{\sqrt{d_0{\K}}}\,
\arctan \frac{\sqrt{\sqrt {\K}-{\bb}_0}}{\sqrt{2d_0}}\,p-
\frac{\sqrt 2\,\alpha_0\sqrt{\sqrt {\K}+{\bb}_0}}{\sqrt{d_0{\K}}}
\,\ln \frac{|p-\sqrt{\frac{2d_0}{\sqrt {\K}+{\bb}_0}}|}
{p+\sqrt{\frac{2d_0}{\sqrt {\K}+{\bb}_0}}}\,,$$
$$0<p^2<\frac{2d_0}{\sqrt {\K}+{\bb}_0}.$$
\vskip 2mm
Type 2. \quad ${\K}>0, \quad \displaystyle{-\frac{16\alpha_0^4+{\K}}
{8\alpha_0^2}}\leq {\bb}_0<-\sqrt {\K} \quad ({\bb}_0^2-{\K}>0)$.

$$t=\frac{2\sqrt 2\,\alpha_0\sqrt{\sqrt {\K}-{\bb}_0}}{\sqrt{d_0{\K}}}\,
\arctan \frac{\sqrt{\sqrt {\K}-{\bb}_0}}{\sqrt{2d_0}}\,p-
\frac{2\sqrt 2\,\alpha_0\sqrt{-(\sqrt {\K}+{\bb}_0)}}{\sqrt{d_0{\K}}}\,
\arctan \frac{\sqrt{-(\sqrt {\K}+{\bb}_0)}}{\sqrt{2d_0}}\,p\,,$$

$$p^2>0.$$
\vskip 2mm
Type 3.\quad ${\K}>0, \quad {\bb}_0>\sqrt {\K} \quad ({\bb}_0^2-{\K}>0)$.

$$t=\frac{\sqrt 2\,\alpha_0\sqrt{{\bb}_0-\sqrt {\K}}}{\sqrt{d_0{\K}}}
\,\ln \frac{|p-\sqrt{\frac{2d_0}{{\bb}_0-\sqrt {\K}}}|}
{p+\sqrt{\frac{2d_0}{{\bb}_0-\sqrt {\K}}}}
-\frac{\sqrt 2\,\alpha_0\sqrt{{\bb}_0+\sqrt {\K}}}{\sqrt{d_0{\K}}}
\,\ln \frac{|p-\sqrt{\frac{2d_0}{{\bb}_0+\sqrt {\K}}}|}
{p+\sqrt{\frac{2d_0}{{\bb}_0+\sqrt {\K}}}}\,,$$

$$p^2 \in \left(0,\,\frac{2d_0}{{\bb}_0+\sqrt {\K}}\right)
\cup \left(\frac{2d_0}{{\bb}_0-\sqrt {\K}},\;\infty\right).$$
\vskip 2mm
Type 4.\quad ${\K}=0, \quad -2\alpha_0^2\leq {\bb}_0<0.$

$$t=\frac{4\alpha_0\,p}{2d_0-{\bb}_0\,p^2}+
\frac{2\sqrt 2\alpha_0}{\sqrt{-{\bb}_0d_0}}
\arctan \frac{\sqrt{-{\bb}_0}}{\sqrt{2d_0}}\,p\,, \quad p^2>0.$$
\vskip 2mm
Type 5.\quad ${\K}=0, \quad {\bb}_0>0.$
$$t=\frac{\sqrt 2\alpha_0}{\sqrt{{\bb}_0d_0}}\,
\ln \frac{p+\sqrt{\frac{2d_0}{{\bb}_0}}}{|p-\sqrt{\frac{2d_0}{{\bb}_0}}|}
+\frac{4\alpha_0\,p}{2d_0-{\bb}_0\,p^2}\,,\quad
p^2\in \left(0,\,\frac{2d_0}{{\bb}_0}\right)\cup
\left(\frac{2d_0}{{\bb}_0},\;\infty\right).$$
\vskip 2mm
Type 6.\quad $-16\alpha_0^2d_0\leq {\K}<0, \quad -2\alpha_0^2
-\displaystyle{\frac{{\K}}{8\alpha_0^2}}\leq {\bb}_0.$

$$t=\frac{2\alpha_0}{\mu}\,\ln\frac{\lambda p^2-\mu p+\nu}
{\lambda p^2+\mu p +\nu}+
\frac{4\alpha_0}{\bar\mu}
\left(\arctan\frac{2\lambda p-\mu}{\sqrt{4\lambda\nu-\mu^2}}
+\arctan\frac{2\lambda p+\mu}{\sqrt{4\lambda\nu+\mu^2}}\right)\,,$$
$$p^2>0.$$
where
$$\lambda=\sqrt{{\bb}_0^2-{\K}}, \quad
\mu=2\sqrt {d_0}\sqrt{\sqrt{{\bb}_0^2-{\K}}+{\bb}_0},\quad
\bar \mu=2\sqrt {d_0}\sqrt{\sqrt{{\bb}_0^2-{\K}}-{\bb}_0}, \quad
\nu=2d_0.$$
\vskip 2mm
\centerline {\it II. The case \quad $a+k^2<0$.}
\vskip 2mm
This case is determined by $d_0<0$, which is equivalent to
$H_0<-3\alpha_0^2$.
\vskip 2mm
Type 7.\quad ${\K}\geq -16\alpha_0^2d_0, \quad -2\alpha_0^2-
\displaystyle{\frac{{\K}}{8\alpha_0^2}}\leq {\bb}_0<-\sqrt {\K} \quad
({\bb}_0^2-{\K}>0).$
$$t=\frac{\sqrt 2\,\alpha_0\sqrt{\sqrt {\K}-{\bb}_0}}{\sqrt{-d_0{\K}}}
\,\ln \frac{|p-\sqrt{\frac{2d_0}{{\bb}_0-\sqrt {\K}}}|}
{p+\sqrt{\frac{2d_0}{{\bb}_0-\sqrt {\K}}}}
-\frac{\sqrt 2\,\alpha_0\sqrt{-({\bb}_0+\sqrt {\K})}}{\sqrt{-d_0{\K}}}
\,\ln \frac{|p-\sqrt{\frac{2d_0}{{\bb}_0+\sqrt {\K}}}|}
{p+\sqrt{\frac{2d_0}{{\bb}_0+\sqrt {\K}}}}\,,$$
$$\frac{2d_0}{{\bb}_0-\sqrt {\K}}<p^2<\frac{2d_0}{{\bb}_0+\sqrt {\K}}.$$
\vskip 2mm
Type 8.\quad ${\K}\geq -16\alpha_0^2d_0, \quad -\sqrt {\K}<{\bb}_0<\sqrt {\K}
\quad ({\bb}_0^2-{\K}<0).$

$$t=\frac{\sqrt 2\,\alpha_0\sqrt{\sqrt {\K}-{\bb}_0}}{\sqrt{-d_0{\K}}}
\,\ln \frac{|p-\sqrt{\frac{2d_0}{{\bb}_0-\sqrt {\K}}}|}
{p+\sqrt{\frac{2d_0}{{\bb}_0-\sqrt {\K}}}}-
\frac{2\sqrt 2\,\alpha_0\sqrt{\sqrt {\K}+{\bb}_0}}{\sqrt{-d_0{\K}}}\,
\arctan \frac{\sqrt{\sqrt {\K}+{\bb}_0}}{\sqrt{-2d_0}}\,p\,,$$
$$p^2\in \left(\frac{2d_0}{{\bb}_0-\sqrt {\K}},\;\infty\right).$$
\vskip 2mm

\centerline {\it Warped product Bochner-K\"ahler metrics with \;
$a+k^2=0 \; \; (d_0=H_0+3\alpha_0^2=0).$}
\vskip 2mm

In this case (5.26) gives that ${\bb}_0^2={\K}$, which implies that
the equations (6.2) and (6.3) coincide.

Considering the positive function $p(s)$ with respect to the
parameter $t$ the equality (6.2) becomes
$$\frac{d^2p}{dt^2}-\frac{4}{p}\left(\frac{dp}{dt}\right)^2
-\frac{{\bb}_0}{2\alpha_0}\,p\,\frac{dp}{dt}=0.\leqno{(6.8)}$$
Taking into account the initial conditions $p(t_0)=1$,\;
$\displaystyle{\frac{dp}{dt}(t_0)=\alpha_0}$ we find the general
solution of (6.8)
$$t=4\alpha_0\int\frac{dp}{p^2\{(4\alpha_0^2-{\bb}_0)p^2+{\bb}_0\}}.
\leqno{(6.9)}$$
From (6.9) we obtain the following functions generating warped
product Bochner-K\"ahler metrics:
\vskip 2mm
Type 9.\quad ${\bb}_0={\K}=0.$
$$p(t)=\frac{1}{\sqrt[3]{1-3\alpha_0(t-t_0)}}\,, \quad
t\in \left(-\infty,\,\frac{1+3\alpha_0t_0}{3\alpha_0}\right).$$
\vskip 2mm
Type 10.\quad ${\bb}_0=-\sqrt {{\K}}<0.$
$$t=\frac{4\alpha_0}{-{\bb}_0p}+\frac{2\alpha_0\sqrt{4\alpha_0^2-{\bb}_0}}
{-{\bb}_0\sqrt{-{\bb}_0}}\,\ln\frac{|\sqrt{4\alpha_0^2-{\bb}_0}\;p-
\sqrt{-{\bb}_0}|} {\sqrt{4\alpha_0^2-{\bb}_0}\;p+\sqrt{-{\bb}_0}}\; ,$$
$$p\in\left(0,\,\frac{\sqrt{-{\bb}_0}}{\sqrt{4\alpha_0^2-{\bb}_0}}\right)
\cup\left(\frac{\sqrt{-{\bb}_0}}{\sqrt{4\alpha_0^2-{\bb}_0}},\; \infty\right).$$
\vskip 2mm
Type 11.\quad ${\bb}_0=\sqrt{{\K}}>0,\quad 4\alpha_0^2-{\bb}_0>0.$
$$t=-\frac{4\alpha_0}{{\bb}_0p}-
\frac{\sqrt{4\alpha_0^2-{\bb}_0}}{\sqrt{{\bb}_0}}
\arctan\frac{\sqrt{4\alpha_0^2-{\bb}_0}}{\sqrt{{\bb}_0}}p\,, \quad p>0.$$
\vskip 2mm
Type 12. \quad ${\bb}_0=\sqrt{\K}>0,\quad 4\alpha_0^2-{\bb}_0=0.$
$$p(t)=\frac{1}{1-\alpha_0(t-t_0)}\,, \quad
t\in\left(-\infty,\,\frac{1+\alpha_0t_0}{\alpha_0}\right).$$
\vskip 2mm
Type 13.\quad ${\bb}_0=\sqrt{{\K}}>0,\quad 4\alpha_0^2-{\bb}_0<0.$
$$t=-\frac{4\alpha_0}{{\bb}_0p}-\frac{\sqrt{{\bb}_0-4\alpha_0^2}}
{2\sqrt {\bb}_0}\,\ln\frac{|\sqrt{{\bb}_0-4\alpha_0^2}\;p-\sqrt {\bb}_0|}
{\sqrt{{\bb}_0-4\alpha_0^2}\;p+\sqrt {\bb}_0}\;,$$
$$p\in\left(0,\,\frac{\sqrt {\bb}_0}{\sqrt{{\bb}_0-4\alpha_0^2}}\right)
\cup \left(\frac{\sqrt {\bb}_0}{\sqrt{{\bb}_0-4\alpha_0^2}},
\;\infty\right).$$

Summarizing the above results and combining with Theorem \ref{T:3.1}
we obtain the main theorem in this section.
\begin{thm}\label{T:6.1}
Any Bochner-K\"ahler manifold whose scalar distribution is a
$B_0$-distribution locally has the structure of a warped
product K\"ahler manifold with metric given by $(6.1)$ and
function $p(t)$ $($or $t(p))$ of type $1. - 13.$
\end{thm}
\vskip 2mm
\begin{rem}
In \cite{GM2} we proved that a K\"ahler metric is of quasi-constant
holomorphic sectional curvatures satisfying the condition $a+k^2>0$
if and only if it locally has the form $\partial\bar\partial f(r^2)$,
$r^2$ being the distance function from the origin in ${\C}^n$.
According to Theorem \ref{T:3.1} the Bochner-K\"ahler metrics of type
$\partial\bar\partial f(r^2)$ (see \cite{TL,B}) can be considered as
warped product Bochner-K\"ahler metrics satisfying the condition
$a+k^2>0$. These metrics are described explicitly in the types 1. - 6.
\end{rem}

\begin{rem}
In \cite{GM3} we proved that a K\"ahler metric is of quasi-constant
holomorphic sectional curvatures satisfying the condition $a+k^2<0$
if and only if it locally has the form $\partial\bar\partial f(-r^2)$,
$-r^2$ being the time-like distance function from the origin in
${\mathbb T}^{n-1}_1$. According to Theorem \ref{T:3.1} the
Bochner-K\"ahler metrics of type $\partial\bar\partial f(-r^2)$
can be considered as warped product Bochner-K\"ahler metrics satisfying
the condition $a+k^2<0$. These metrics are described explicitly in the
types 7. - 8.
\end{rem}

Finally we shall describe the complete warped product Bochner-K\"ahler
metrics.

According to \cite{BO} a warped product K\"ahler metric
(6.1) is complete if and only if the base $Q_0$ is a complete
$\alpha_0$-Sasakian space form and the function $f(t)$ is defined in
the interval $I={\R}$. Investigating the functions of type 1. - 13. we
obtain four families of complete metrics:
\vskip 2mm
{\it The family \quad $({\K}=0,\quad {\bb}_0>0, \quad d_0>0)$.}
\vskip 2mm
Given a complete $\alpha_0$-Sasakian space form $Q_0$ satisfying the
condition $d_0=H_0+3\alpha_0^2>0$ (cf \cite{T}) and a constant
${\bb}_0>0$. Then the function (Type 5.)
$$u=\frac{1}{t}=\left\{\frac{\sqrt 2\alpha_0}{\sqrt{{\bb}_0d_0}}\,
\ln \frac{p+\sqrt{\frac{2d_0}{{\bb}_0}}}{|p-\sqrt{\frac{2d_0}{{\bb}_0}}|}
+\frac{4\alpha_0\,p}{2d_0-{\bb}_0\,p^2}\right\}^{-1},\quad p>0$$
determines a one-parameter family (with respect to ${\bb}_0$)
of complete metrics (cf \cite{B}).
\vskip 2mm
{\it The family \quad $({\K}\geq -16\alpha_0^2d_0, \quad -2\alpha_0^2-
\displaystyle{\frac{{\K}}{8\alpha_0^2}}\leq {\bb}_0<-\sqrt {\K}, \quad
d_0<0)$.}
\vskip 2mm
Given a complete $\alpha_0$-Sasakian space form $Q_0$ satisfying the
condition $d_0=H_0+3\alpha_0^2<0$ (cf \cite{T}) and two constants
${\K}\geq -16\alpha_0^2d_0$ and $-2\alpha_0^2-
\displaystyle{\frac{{\K}}{8\alpha_0^2}}\leq {\bb}_0<-\sqrt {\K}$.
Then the function  (Type 7.)
$$t=\frac{\sqrt 2\,\alpha_0\sqrt{\sqrt {\K}-{\bb}_0}}{\sqrt{-d_0{\K}}}
\,\ln \frac{|p-\sqrt{\frac{2d_0}{{\bb}_0-\sqrt {\K}}}|}
{p+\sqrt{\frac{2d_0}{{\bb}_0-\sqrt {\K}}}}
-\frac{\sqrt 2\,\alpha_0\sqrt{-({\bb}_0+\sqrt {\K})}}{\sqrt{-d_0{\K}}}
\,\ln \frac{|p-\sqrt{\frac{2d_0}{{\bb}_0+\sqrt {\K}}}|}
{p+\sqrt{\frac{2d_0}{{\bb}_0+\sqrt {\K}}}}\,,$$
$$\frac{2d_0}{{\bb}_0-\sqrt {\K}}<p^2<\frac{2d_0}{{\bb}_0+\sqrt {\K}}$$
determines a two-parameter family (with respect to ${\bb}_0$ and
${\K}$) of complete metrics.
\vskip 2mm
{\it The family \quad $({\bb}_0=-\sqrt{\K}<0,\quad d_0=0)$.}
\vskip 2mm
Given a complete $\alpha_0$-Sasakian space form $Q_0$ satisfying the
condition $d_0=H_0+3\alpha_0^2=0$ (cf \cite{T}) and a constant
${\bb}_0<0$. Then the function (Type 10.)
$$t=\frac{4\alpha_0}{-{\bb}_0p}+\frac{2\alpha_0\sqrt{4\alpha_0^2-{\bb}_0}}
{-{\bb}_0\sqrt{-{\bb}_0}}\ln\frac{|\sqrt{4\alpha_0^2-{\bb}_0}\;p-
\sqrt{-{\bb}_0}|}{\sqrt{4\alpha_0^2-{\bb}_0}\;p+\sqrt{-{\bb}_0}}\,, \quad
p\in\left(0,\,\frac{\sqrt{-{\bb}_0}}{\sqrt{4\alpha_0^2-{\bb}_0}}\right)$$
determines a one-parameter family (with respect to ${\bb}_0$) of
complete metrics.
\vskip 2mm
{\it The family \quad $(4\alpha_0^2<{\bb}_0=\sqrt{\K},\quad d_0=0)$.}
\vskip 2mm
Given a complete $\alpha_0$-Sasakian space form $Q_0$ satisfying the
condition $d_0=H_0+3\alpha_0^2=0$ and a constant
${\bb}_0>4\alpha_0^2$. Then the function (Type 13.)
$$t=-\frac{4\alpha_0}{{\bb}_0p}-\frac{\sqrt{{\bb}_0-4\alpha_0^2}}
{2\sqrt {\bb}_0}\ln\frac{|\sqrt{{\bb}_0-4\alpha_0^2}\;p-\sqrt {\bb}_0|}
{\sqrt{{\bb}_0-4\alpha_0^2}\;p+\sqrt {\bb}_0}\,,\quad
p\in\left(0,\,\frac{\sqrt {\bb}_0}{\sqrt{{\bb}_0-4\alpha_0^2}}\right)$$
determines a one-parameter family (with respect to ${\bb}_0$) of
complete metrics.
\vskip 2mm
\section*{Notes on the geometry of K\"ahler manifolds with $J$-invariant
distributions of codimension two}

We note that every K\"ahler manifold $(M,g,J,D) \;
(\dim \, M = 2n \geq 4)$
with $J$-invariant distribution $D$ carries the tensor
$$P=\frac{2}{(n+1)(n+2)}\,\pi - \frac{4}{n+2}\, \Phi + \Psi,$$
which is the unique (up to a factor) invariant tensor with
zero Ricci trace.

We can draw a parallel between the K\"ahler manifolds
$(M,g,J)$ whose structural group is $U(n)$ and the K\"ahler
manifolds $(M,g,J,D)$ whose structural group is $U(n-1)\times
U(1)$. We compare a class of K\"ahler manifolds $(M,g,J,D)$ with
any of the basic classes of K\"ahler manifolds $(M,g,J)$. The
correspondence between the curvature identities, which
characterize these classes is given as follows
$$\begin{array}{rcl}
(M,g,J) & \longleftrightarrow & (M,g,J,D)\\
[2mm]\displaystyle{R=\frac{\tau}{n(n+1)}\,\pi} &
\longleftrightarrow & R=a\pi +b\Phi+c\Psi, \\
[4mm]
B(R)=0 & \longleftrightarrow & B(R)=cP,\\
[2mm] \displaystyle{\rho = \frac{\tau}{2n}\,g} & \longleftrightarrow
& \displaystyle{\rho = \frac{\tau - 2\sigma}{2(n-1)}\,g +
\frac{2n\sigma -
\tau}{2(n-1)}\,(\eta \otimes \eta + \tilde\eta \otimes \tilde\eta)}.\\
\end{array}$$

The next natural proposition is valid.
\vskip 2mm
{\it A K\"ahler manifold $(M,g,J,D) \; (\dim \,M=2n\geq4)$ is
of quasi-constant holomorphic sectional curvatures if and only
if the following curvature identities hold good}
$$\begin{array}{ll}
\vspace{2mm}
(i) & B(R)=cP;\\
\vspace{2mm}
(ii) & \displaystyle{\rho = \frac{\tau - 2\sigma}{2(n-1)}\,g +
\frac{2n\sigma - \tau}{2(n-1)}\,(\eta \otimes \eta
+ \tilde\eta \otimes \tilde\eta)}.
\end{array}$$

Let $(M,g,J)$ be a K\"ahler manifold with $d\tau\neq 0$
and $\eta=\displaystyle{\frac{d\,\tau}{\Vert d\,\tau \Vert}}$.

Suppose that
$$ 1)\; (grad\;\tau)^{1,0} \; \; \rm{is \; a \; holomorphic \;vector
\; field;} \quad
2) \; \nabla_{\lambda}\,\eta_{\bar\mu}=\frac{k}{2} \, g_{\lambda\bar\mu}.
\leqno{(*)}$$
Under the conditions (*) any biconformal transformation of the
structure $(g,\eta)$ again gives a K\"ahler structure
(cf \cite{GM2}).

In K\"ahler geometry it seems that the following statement is true.
\vskip 2mm
{\it Let $(M,g,J)$ be a K\"ahler manifold whose scalar distribution
satisfies the conditions $(*)$. Then the tensor $B(R)-cP$ is a
biconformal invariant.}


\begin{thebibliography}{99}
\bibitem{BO}
Bishop, R.; O'Neil, B. {\it Manifolds of negative curvature,}
Trans. Amer. Math. Soc., {\bf 145} (1969), 1-49.
\bibitem{BP}
Boju, V.; Popesku, M. {\it Espaces \`a courbure quasi-constante,}
J. Diff. Geom., {\bf 13} (1978), 373-383.
\bibitem{B}
Bryant, R. {\it Bochner-K\"ahler metrics,} J. Amer. Math. Soc.,
{\bf 14} (2001), 623-715.
\bibitem{GH}
Gray, A.; Hervella, L. {\it The sixteen classes of almost
Hermitian manifolds and their linear invariants}, Ann. di Mat.
Pura ed Appl., {\bf 123} (1980), 35-58.
\bibitem{GM1}
Ganchev, G.; Mihova, V. {\it Riemannian manifolds of
quasi-constant sectional curvature}, J. reine und angew. Math.,
{\bf 522} (2000), 119-141.
\bibitem{GM2}
Ganchev, G.; Mihova, V. {\it K\"ahler manifolds of quasi-constant
holomorphic sectional curvatures}, ArXiv: math.DG/0505671, to appear.
\bibitem{GM3}
Ganchev, G.; Mihova, V. {\it K\"ahler metrics generated by
functions of the time-like distance in the flat K\"ahler-Lorentz
space}, ArXiv: math.DG/0510468, to appear.
\bibitem{JV}
Janssens, D.; Vanhecke, L. {\it Almost contact structures and
curvature tensors}, Kodai Math. J., {\bf 4} (1981), 1-27.
\bibitem{K}
Kamishima, Y. {\it Uniformization of K\"ahler manifolds with
vanishing Bochner tensor}, Acta Math., {\bf 172} (1994), 299-308.
\bibitem{KN}
Kobayasi, S.; Nomizu, K. {\it Foundations of Differential
Geometry, II,} Interscience Publishers, New Yourk, 1969.
\bibitem{O}
Ogiue, K. {\it On almost contact manifolds admitting axiom
of planes or axiom of free mobility,}
Kodai Math. Sem. Rep., {\bf 16} (1964), 223-232.
\bibitem{Ok}
Okumura, M. {\it On infinitesimal conformal and projective
transformations of normal contact spaces,}
T\^{o}hoku Math. J., {\bf 14} (1962), 398-412.
\bibitem{TL}
Tachibana, S.; Liu, R.C. {\it Notes on K\"ahlerian metrics with
vanishing Bochner curvature tensor}, Kodai Math. Sem. Rep.,
{\bf 22} (1970), 313-321.
\bibitem{T}
Tanno, S. {\it Sasakian manifolds with constant
$\varphi$-holomorphic sectional curvature,} T\^ohoku Math. J.,
{\bf 21} (1969), 501-507.
\bibitem{T1}
Tashiro, Y. {\it On contact structures on hypersurfaces in almost
complex manifolds I,} T\^ohoku Math. J., {\bf 15} (1963), 62-79.
\bibitem{T2}
Tashiro, Y. {\it On contact structures on hypersurfaces in almost
complex manifolds II,} T\^ohoku Math. J., {\bf 15} (1963), 167-175.
\bibitem{TV}
Tricerri, F.; Vanhecke, L. {\it Curvature tensors on almost
Hermitian manifolds}, Trans. Amer. Math. Soc., {\bf 267} (1981),
365-398.
\end{thebibliography}
\end{document}